\documentclass{amsart}
\usepackage{amssymb,latexsym}
\usepackage{amsmath,amsfonts,amsthm}
\usepackage[dvips]{graphicx}
\newtheorem{theorem}{Theorem}[section]

\newtheorem{lemma}[theorem]{Lemma}
\newtheorem{remark}[theorem]{Remark}

\newtheorem{definition}[theorem]{Definition}
\newtheorem{cor}[theorem]{Corollary}
\def\co{\colon\thinspace}

\title{Kodaira dimension and symplectic sums}
\author{Michael Usher}\address{Department of Mathematics\\ University of Georgia\\ Athens, GA  30606}\email{usher@math.uga.edu}
\begin{document}
\begin{abstract}
Modulo trivial exceptions, we show that  symplectic sums  of symplectic
$4$-manifolds along surfaces of positive genus are never rational or
ruled, and we enumerate each case in which they have Kodaira
dimension zero (\emph{i.e.}, are blowups of symplectic $4$-manifolds
with torsion canonical class).  In particular, a symplectic
four-manifold of Kodaira dimension zero arises by such a surgery
only if it is diffeomorphic to a blowup either of the $K3$ surface,
the Enriques surface, or a member of a particular family of
$T^2$-bundles over $T^2$ each having $b_1=2$.
\end{abstract}
\maketitle
\section{Introduction}

Our understanding of the diversity of the world of symplectic
four-manifolds has been greatly enriched by the introduction in
\cite{G} and \cite{MW} of the \emph{symplectic sum}. Given
symplectic four-manifolds $(X_1,\omega_1),(X_2,\omega_2)$ containing
embedded, two-dimensional symplectic submanifolds $F_1\subset X_1$,
$F_2\subset X_2$ of equal area and genus and an
orientation-reversing isomorphism $\Phi\co N_{X_1}F_1\to N_{X_2}F_2$
of their normal bundles (which of course exists if and only if $F_1$
and $F_2$ have opposite self-intersection), the symplectic sum
operation provides a natural isotopy class of symplectic structures
on the normal connect sum \[ Z=X_1\#_{F_1=F_2}X_2=(X_1\setminus
\nu_1)\cup_{\partial\nu_1\sim_{\Phi}\partial\nu_2}(X_2\setminus
\nu_2),\] where the $\nu_i$ are tubular neighborhoods of $F_i$ and
we use the restriction of $\Phi$ to the unit normal circle bundles of the
$F_i$ to glue the boundaries of the manifolds $X_i\setminus \nu_i$.
Using the symplectic sum along surfaces of positive genus, various
authors over the years have constructed symplectic four-manifolds
satisfying an impressive array of properties; see for instance
Theorem 6.2 of \cite{G}, which for any finitely presented group $G$
gives a number $r(G)$ such that whenever $a+b\equiv 0 \mod 12$ and
$0\leq a\leq 2(b-r(G))$ there is a symplectic $4$-manifold
$M_{a,b,G}$ with $\pi_1(M_{a,b,G})=G$, $c_{1}^{2}(M_{a,b,G})=a$, and
$c_2(M_{a,b,G})=b$.  While Gompf's examples were distinguished by
their classical topological invariants, the symplectic sum also
gives rise to infinite families of mutually homeomorphic but
nondiffeomorphic symplectic four-manifolds, since if $K$ is a
fibered knot the operation of knot surgery with $K$ \cite{FS}
amounts to a symplectic sum.

The purpose of this note is to show that, notwithstanding the
diversity of symplectic four-manifolds that can be constructed via
symplectic sum, there are  significant topological restrictions on
the manifolds that can be obtained in this way.  Our results may
perhaps best be understood in terms of the notion of (symplectic)
\emph{Kodaira dimension}, a notion which dates from \cite{MSb} and
is discussed in some detail in \cite{Li}.  We shall recall the
definition of Kodaira dimension below. First, recall that a
symplectic four-manifold is called \emph{minimal} if it does not
contain any embedded symplectic spheres of square $-1$, and that if
$(X,\omega)$ is any symplectic four-manifold one may obtain a
minimal symplectic four-manifold $(X',\omega')$ by blowing down a
maximal disjoint collection of symplectic $(-1)$-spheres in $X$;
$(X',\omega')$ is then called a minimal model for $(X,\omega)$.

\begin{definition}  Let $(X,\omega)$ be a symplectic four-manifold
with minimal model $(X',\omega')$, and let $\kappa_{X'}\in
H^2(X';\mathbb{Z})$ denote the canonical class of $(X',\omega')$.
Then the Kodaira dimension of $(X,\omega)$ is \[
\kappa(X,\omega)=\left\{\begin{array}{rl}-\infty& \mbox{if
}\kappa_{X'}\cdot[\omega']<0 \mbox{ or } \kappa_{X'}^{2}<0 \\
0 & \mbox{if } \kappa_{X'}\cdot[\omega']=\kappa_{X'}^{2}=0\\ 1 &
\mbox{if }\kappa_{X'}\cdot[\omega']>0\mbox{ and
}\kappa_{X'}^{2}=0\\2& \mbox{if }\kappa_{X'}\cdot[\omega']>0\mbox{
and }\kappa_{X'}^{2}>0\end{array}\right.
\]
\end{definition}

In Section 2 of \cite{Li} and references therein it is shown that
$\kappa(X,\omega)$ is well-defined for any symplectic four-manifold
(in particular it is independent of the choice of minimal model, and one of the four possibilities listed above always holds);
coincides with the classical notion of Kodaira dimension in cases
when $X$ happens to admit the structure of a complex surface;  is
equal to $-\infty$ if and only if $X$ is a rational or ruled
surface;\footnote{Here and below we adopt the convention that a
ruled surface is a symplectic manifold obtained by (possibly)
blowing up the total space of an $S^2$-bundle over some Riemann
surface. When we want to assume that no blowups have been carried
out, we shall refer instead to an ``$S^2$-bundle.''} and is equal to
zero if and only if the canonical class of the minimal model $X'$ is
torsion. Moreover Theorem 2.6 of \cite{Li} shows that the Kodaira
dimension $\kappa(M,\omega)$ depends only on the diffeomorphism type
of $M$, and not on the symplectic form.

A common feature of the numerous interesting new symplectic
four-manifolds that the symplectic sum operation has provided to us
is that they have always had Kodaira dimension 1 or 2.  For
instance, knot surgery on the K3 surface with a nontrivial fibered
knot as in \cite{FS} always yields a symplectic four-manifold of
Kodaira dimension 1 (even though the result is homeomorphic to the
K3 surface, which has Kodaira dimension zero), and (at least aside
from some very trivial cases) Gompf's manifolds $M_{a,b,G}$ have
Kodaira dimension 1 if $a=0$ and 2 if $a>0$.  Our main theorem below
will demonstrate that this is not a coincidence.  To state it, we
make the following  definitions.

\begin{definition} Let $(X_1,\omega_1)$, $(X_2,\omega_2)$ be 
symplectic four-manifolds with $F_i\subset X_i$ embedded symplectic
submanifolds of equal area and genus and opposite square. \begin{enumerate}\item The
symplectic sum $X_1\#_{F_1=F_2}X_2$ is called \emph{smoothly
trivial} if, for some $i\in\{1,2\}$, $X_i$ the total space of an
$S^2$-bundle of which $F_i$ is a section.  Otherwise, the symplectic
sum is called smoothly nontrivial.\item The symplectic sum $X_1\#_{F_1=F_2}X_2$ is said to be \emph{of blowup type} if, for some $i\in\{1,2\}$, the pair $(X_i,F_i)$ may be obtained from a pair $(E,F)$ consisting of the total space $E$ of an $S^2$-bundle of which $F$ is a section by a sequence of zero or more blowups at points not lying on $F$.\end{enumerate}

\end{definition}

\begin{remark}\label{merid}  It is stated without proof in \cite{G} that if
$X_1\#_{F_1=F_2}X_2$ is a smoothly trivial symplectic sum, say with
$(X_2,F_2)$ consisting of an $S^2$-bundle and a section, then the
sum $X_1\#_{F_1=F_2}X_2$ is diffeomorphic to $X_1$.  It is not
difficult to prove this: simply note that $X_2\setminus \nu_2$ will
be diffeomorphic to a neighborhood of $F_1$ in $X_1$, so that for at
least one choice of the gluing map $\Phi|_{\partial\nu_1}$ the sum
will be diffeomorphic to $X_1$; further, the gluing maps
$\Phi|_{\partial\nu_1}\co
\partial\nu_1\to\partial\nu_2$ that we are allowed to use in forming the
symplectic sum are precisely the restrictions of orientation
reversing bundle isomorphisms, and so any two of them differ by
precomposing with an orientation preserving diffeomorphism of
$\partial\nu_1$ which extends over $\nu_1$, implying therefore that
the diffeomorphism type of $X_1\#_{F_1=F_2}X_2$ is independent of
the gluing map and so is $X_1$ in any event.  Of course, this
argument is dependent on the fact that the gluing map is  required
to preserve the fibers of the normal circle bundles, an issue which
seems to have caused a certain amount of confusion in the
literature, where one occasionally finds mistaken claims that
symplectic sums with $S^2$-bundles along sections sometimes change
the diffeomorphism type.

At any rate, the above fact justifies our use of the term ``smoothly
trivial'' to describe such symplectic sums, and by Theorem 2.6 of
\cite{Li} implies that performing a smoothly trivial symplectic sum
leaves the Kodaira dimension unchanged.  Incidentally, while these sums are trivial from a smooth standpoint, they generally do alter the symplectic structure in a manner equivalent to the ``inflation'' technique of \cite{LM}, a fact which is exploited in \cite{LU}.

\end{remark}

\begin{definition} If $(X_1,\omega_1)$, $(X_2,\omega_2)$ are
symplectic four-manifolds with $F_i\subset X_i$ embedded symplectic
submanifolds of equal area and genus and opposite square, the
symplectic sum $X_1\#_{F_1=F_2}X_2$ is called \emph{relatively
minimal} if for each $i=1,2$, there are no embedded symplectic
spheres of square $-1$ in $X_i\setminus F_i$.
\end{definition}

\begin{remark}  If a symplectic sum $Z=X_1\#_{F_1=F_2}X_2$ is not
relatively minimal, then if we blow down maximal disjoint collections of spheres of
square $-1$ in $X_1\setminus F_1$ and $X_2\setminus F_2$ to obtain
symplectic manifolds $X'_1,X'_2$, then the $F_i$ survive in the
blowdowns and the areas, genera, and self-intersections of the $F_i$
are left unchanged.  Hence we may form a symplectic sum
$Z'=X'_1\#_{F_1=F_2}X'_2$, and $Z$ may be recovered from $Z'$ by a
sequence of blowups.  $X'_1\#_{F_1=F_2}X'_2$ will be smoothly trivial if and only if $X_1\#_{F_1=F_2}X_2$ is of blowup type.  Thus any symplectic $4$-manifold which arises
as a symplectic sum which is not of blowup type is a blowup of a
symplectic $4$-manifold which arises as a smoothly nontrivial,
relatively minimal symplectic sum.  Moreover, symplectic sums which are of blowup type are diffeomorphic to blowups of one of their summands.\end{remark}

\begin{remark} In this language, the main result of \cite{U} may be
rephrased as stating that any symplectic $4$-manifold arising as a
smoothly nontrivial, relatively minimal symplectic sum along
surfaces of positive genus is minimal.
\end{remark}

Our main result is the following:

\begin{theorem}\label{main}
Suppose that $(X_1,\omega_1)$, $(X_2,\omega_2)$, $F_1$, and $F_2$
are such that the symplectic sum $Z=X_1\#_{F_1=F_2} X_2$ is smoothly
nontrivial and relatively minimal and the genus of the $F_i$ is
positive. Then:\begin{itemize}
\item[(a)] $Z$ does not have Kodaira dimension $-\infty$.
\item[(b)] If $Z$ has Kodaira dimension $0$, then the diffeomorphism
types of $X_1$, $X_2$, and $Z$ are given by one of the rows in the
following table, where the notation $M(A,B;\vec{v})$ denotes a
$T^2$-bundle over $T^2$ as described below or in
\cite{SF}.  Moreover, each entry in the third column of this table can in fact be constructed as a smoothly nontrivial, relatively minimal symplectic sum along a torus.\end{itemize}\begin{center}
\begin{tabular}{|c|c|p{3 in}|}\hline
$X_1$&$X_2$&possible diffeomorphism types of $X_1\#_{F_1=F_2}X_2$\\
\hline \rule{0pt}{.4 cm}
$\mathbb{C}P^2\#(18-k)\overline{\mathbb{C}P^2}$ &
$\mathbb{C}P^2\#k\overline{\mathbb{C}P^2}$ & $K3$ surface \\
\hline \rule{0pt}{.4 cm}$S^2\times S^2$ &
$\mathbb{C}P^2\#17\overline{\mathbb{C}P^2}$ & $K3$ surface
\\ \hline \rule{0pt}{.4 cm}
$\mathbb{C}P^2\#(9-k)\overline{\mathbb{C}P^2}$ & $(S^2\times
T^2)\#k\overline{\mathbb{C}P^2}$ & Enriques surface
\\ \hline \rule{0pt}{.4 cm} $S^2\times S^2$ & $(S^2\times
T^2)\#8\overline{\mathbb{C}P^2}$ & Enriques surface \\ \hline
\rule{0pt}{.4 cm} $\mathbb{C}P^2\#9\overline{\mathbb{C}P^2}$ &
$S^2\tilde{\times}T^2$ & Enriques surface \\ \hline $S^2\times T^2$
& $S^2\times T^2$ &
$M\left(I,\left(\begin{array}{cc}-1&z\\0&-1\end{array}\right);\left(\begin{array}{c}0\\0\end{array}\right)\right)$,
$M\left(-I,\left(\begin{array}{cc}1&2y\\0&1\end{array}\right);\left(\begin{array}{c}0\\0\end{array}\right)\right)$,
$M\left(I,\left(\begin{array}{cc}-1&2y\\0&-1\end{array}\right);\left(\begin{array}{c}0\\1\end{array}\right)\right)$,
$M\left(-I,\left(\begin{array}{cc}1&2y\\0&1\end{array}\right);\left(\begin{array}{c}0\\1\end{array}\right)\right)$
$(y,z\in\mathbb{Z})$\\  \hline $S^2\times T^2$ & $S^2\tilde{\times}
T^2$ &
$M\left(-I,\left(\begin{array}{cc}1&2y+1\\0&1\end{array}\right);\left(\begin{array}{c}0\\0\end{array}\right)\right)$,
$M\left(-I,\left(\begin{array}{cc}1&2y+1\\0&1\end{array}\right);\left(\begin{array}{c}0\\1\end{array}\right)\right)$
$(y\in\mathbb{Z})$\\ \hline  $S^2\tilde{\times} T^2$ &
$S^2\tilde{\times} T^2$ &
$M\left(I,\left(\begin{array}{cc}-1&z\\0&-1\end{array}\right);\left(\begin{array}{c}1\\0\end{array}\right)\right)$,
$M\left(-I,\left(\begin{array}{cc}1&2y\\0&1\end{array}\right);\left(\begin{array}{c}1\\0\end{array}\right)\right)$
$M\left(I,\left(\begin{array}{cc}-1&2y+1\\0&-1\end{array}\right);\left(\begin{array}{c}0\\1\end{array}\right)\right)$
$(y,z\in\mathbb{Z})$\\ \hline
 \end{tabular}\end{center}
\end{theorem}

Here $S^2\tilde{\times}T^2$ refers to the smoothly nontrivial
$S^2$-bundle over the torus. The $K3$ surface is, of course, the
smooth four-manifold underlying a quartic hypersurface of
$\mathbb{C}P^3$, while the Enriques surface can be realized either
as the quotient of a $K3$ surface by a free involution which is
holomorphic for an appropriate complex structure, or, equivalently,
as the result of performing logarithmic transformations of
multiplicity 2 on two fibers of the rational elliptic surface
$E(1)$; pp. 590--599 of \cite{GH} provide a concise description of
the complex geometry of these manifolds. The notation
$M(A,B;\vec{v})$ ($A,B\in SL(2;\mathbb{Z}),\vec{v}\in\mathbb{Z}^2, AB=BA$)
above is as in \cite{SF}: identifying
$T^2=\mathbb{R}^2/\mathbb{Z}^2$ and writing the coset of
$(x,y)\in\mathbb{R}^2$ as
$\left[\begin{array}{c}x\\y\end{array}\right]$, we define \[
M\left(A,B;\left(\begin{array}{c}0\\0\end{array}\right)\right)=\frac{\mathbb{R}^2\times
T^2}{\begin{array}{rl}\left(s+1,t,\left[\begin{array}{c}x\\y\end{array}\right]\right)&\sim
\left(s,t,\left[A\left(\begin{array}{c}x\\y\end{array}\right)\right]\right),\\\left(s,t+1,\left[\begin{array}{c}x\\y\end{array}\right]\right)&\sim
\left(s,t,\left[B\left(\begin{array}{c}x\\y\end{array}\right)\right]\right)\end{array}},\]
and let $M(A,B;(m,n))$ be obtained from $M(A,B;(0,0))$ by removing a
neighborhood of a torus fiber of the projection $[s,t,x,y]\mapsto
[s,t]$ and regluing it by the map \begin{align*} \partial D^2\times
T^2&\to \partial D^2\times T^2\\ (e^{2\pi i\theta},[x,y])&\mapsto
(e^{2\pi i\theta},[x+m\theta],[y+n\theta]).\end{align*}

Since manifolds obtained by general symplectic sums are blowups of
those obtained by relatively minimal ones, and since Kodaira
dimension is preserved under blowups, Theorem \ref{main} has the
following corollary:

\begin{cor}  No symplectic $4$-manifold of Kodaira dimension $-\infty$ arises as a
 symplectic sum along surfaces of positive genus which is not of blowup type,
and the only symplectic $4$-manifolds of Kodaira dimension $0$ which
arise as symplectic sums which are not of blowup type along surfaces of
positive genus are blowups of symplectic manifolds diffeomorphic to
either the $K3$ surface, the Enriques surface, or a $T^2$-bundle
over $T^2$ having the form $M\left(I,
\left(\begin{array}{cc}-1&z\\0&-1\end{array}\right);\vec{v}\right)$
or $M\left(-I,
\left(\begin{array}{cc}1&z\\0&1\end{array}\right);\vec{v}\right)$,
where $z\in \mathbb{Z}$ and $\vec{v}\in\{(0,0),(0,1),(1,0)\}$.
\end{cor}

Symplectic $4$-manifolds of Kodaira dimension $-\infty$ are of
course well understood since they are all rational or ruled.  In
light of Gompf's manifolds $M_{a,b,G}$ mentioned above, the classes
of symplectic $4$-manifolds of Kodaira dimension either $1$ or $2$
are each at least as complicated as the class of all finitely
presented groups, which  have been known to be unclassifiable since
the 1950's.  The classification of symplectic $4$-manifolds of
Kodaira dimension zero, meanwhile, is a very interesting open
problem.  The only presently-known such manifolds are blowups of the
$K3$ surface, the Enriques surface, or an orientable $T^2$-bundle
over $T^2$ (all of the latter were proven to admit symplectic
structures in \cite{Geiges}), and a recent striking result
independently proven in \cite{B} and \cite{Li2} implies that any
minimal symplectic $4$-manifold of Kodaira dimension zero
necessarily has the same rational homology as one of the known
minimal examples. Theorem \ref{main} thus shows that (up to
diffeomorphism) if any new symplectic $4$-manifolds of Kodaira
dimension zero do exist, then they cannot be found by performing
symplectic sums along surfaces of positive genus.  Note that since
blowing down a $(-1)$-sphere amounts to performing a symplectic sum
with $\mathbb{C}P^2$ along a line, any symplectic $4$-manifold
arises (in a somewhat trivial fashion) as a smoothly nontrivial
symplectic sum along spheres.  It would be interesting to know
whether or not new symplectic $4$-manifolds of Kodaira dimension
zero can arise by symplectic sum with $\mathbb{C}P^2$ along a
quadric (\emph{i.e.}, by blowing down a $(-4)$-sphere), or more
generally whether or not new symplectic $4$-manifolds of Kodaira
dimension zero can arise via the rational blowdown procedure of
\cite{FS2}.

Theorem \ref{main} pins down the diffeomorphism type of the
symplectic sums in question; it would be somewhat preferable to have
a result specifying the symplectic deformation type (though it's
worth noting that this may be a moot point, since any two symplectic
structures on any of the manifolds in the last column of the table
in Theorem \ref{main} have the same canonical class, and there are
not currently any known examples of smooth four-manifolds admitting
deformation inequivalent symplectic structures that have the same
canonical class).  In the case that the symplectic sum is a
$T^2$-bundle over $T^2$, our proof can be seen to imply that the
symplectic form will be positive on the fibers, and so will be
deformation equivalent to a form obtained by the Thurston trick.
When the sum is diffeomorphic to the $K3$ surface or the Enriques
surface, the situation is perhaps more subtle, since in our proofs
we use, among other things, the fact that any orientation-preserving
diffeomorphism of the boundary of the complement of a tubular
neighborhood of a fiber in the rational elliptic surface $E(1)$
extends to the whole complement, and so since such a diffeomorphism
cannot always be made symplectic we lose control of the symplectic
form.

We should mention that, as in Remark \ref{merid}, although the
conclusion of Theorem \ref{main} only concerns diffeomorphism types
the result is sensitive to the fact that the sum operation is
symplectic and so the gluing map $\Phi\co \partial
\nu_1\to\partial\nu_2$ is required to be the restriction of an
anti-isomorphism of complex line bundles (rather than just an
arbitrary diffeomorphism).  Of course, if $\Phi$ were replaced by an
arbitrary diffeomorphism, the sum might not admit a symplectic
structure; for instance for $n\geq 2$ a degree-zero logarithmic
transformation of a fiber of the elliptic surface $E(n)$ (which results from
an appropriate gluing of the complement of the symplectic torus
fiber of $E(n)$ to the complement of a symplectic torus in
$T^2\times S^2$) is diffeomorphic to the non-symplectic manifold
$(2n-1)\mathbb{C}P^2\#(10n-1)\overline{\mathbb{C}P^2}$ (Theorem 6.1
of \cite{G1}).  In other cases, though, a more general gluing map
gives rise to a manifold that happens to admit a symplectic
structure, which may even have Kodaira dimension zero: for instance,
any of the $T^2$-bundles over $T^2$ given in our notation above as
$M(I,I;(m,n))$ may be obtained by an appropriate gluing of the
complement of a symplectic torus in $T^4$ to the complement of a
symplectic torus in $T^2\times S^2$; these all admit symplectic
structures by \cite{Geiges} (or indeed by earlier work; the case
$(m,n)=(0,1)$ is the Kodaira-Thurston manifold) with trivial
canonical class and so Kodaira dimension zero.  However, for
$(m,n)\neq (0,0)$ $b_1(M(I,I;(m,n)))=3$, whereas all of the
$T^2$-bundles over $T^2$ appearing in Theorem \ref{main} have
$b_1=2$, so none of the manifolds $M(I,I;(m,n))$ occur as smoothly
nontrivial symplectic sums.

In the next section, we prove part (a) of Theorem \ref{main}, which
might  best be seen as a corollary of the arguments of \cite{U}. The
rest of the paper consists of the proof of part (b), which relies on
a result describing the behavior of the canonical class of a
symplectic sum in order to reduce the result to a small number of
cases which can be dispensed with in turn.  A useful tool in two of
the cases is Lemma \ref{isotdiff}, which states that if $F$ is a
symplectic torus  in
$\mathbb{C}P^2\#9\overline{\mathbb{C}P^2}$ Poincar\'e dual to the
first Chern class, then there is a diffeomorphism of pairs taking
$(\mathbb{C}P^2\#9\overline{\mathbb{C}P^2},F)$ to the pair
consisting of a rational elliptic surface together with one of its
regular fibers.  This should be contrasted with various non-isotopy
results for symplectic tori Poincar\'e dual to multiples of the
first Chern class, \emph{e.g.}, in \cite{FS3}.
\subsection*{Acknowledgement} I am grateful for the referee's corrections and suggestions.
This work was partially supported by an NSF Postdoctoral Fellowship, and was carried out while I was at Princeton University.

\section{Rational and ruled surfaces}

Part (a) of Theorem \ref{main} is an easy consequence of two results
proven in \cite{U}.  On the one hand, from Section 2 of that paper,
we find (based largely on \cite{IP}), that if $Z=X_1\#_{F_1=F_2}
X_2$  is a symplectic sum, and if $A\in H_2(Z;\mathbb{Z})$ is any
nonzero homology class having a nontrivial Gromov--Witten invariant
counting genus-zero holomorphic curves, then there are classes
$A_1\in H_2(X_1;\mathbb{Z})$, $A_2\in H_2(X_2;\mathbb{Z})$, each
represented by a union of genus zero stable maps which are pseudoholomorphic
for some almost complex structure $J_i$ on $X_i$ ($i=1,2$) which
preserves $TF_i$ such that
\begin{equation}\label{minusone}
\langle\kappa_{X_1}+PD[F_1],A_1\rangle+\langle\kappa_{X_2}+PD[F_2],A_2\rangle=\langle \kappa_Z,A\rangle<0,\end{equation}
where the final inequality results from the fact that the expected dimension of genus-zero pseudoholomorphic representatives of $A$, namely $-1-\langle\kappa_Z,A\rangle$, must be nonnegative in order for the corresponding Gromov--Witten invariant to be nonvanishing.

On the other hand, Proposition 3.9 of \cite{U}, translated into the
terminology of the introduction, states that if the symplectic sum
$Z=X_1\#_{F_1=F_2} X_2$ is smoothly nontrivial and relatively
minimal, and if the $F_i$ have positive genus, then the $F_i$ are
``rationally $K$-nef,'' which is to say, for $i=1,2$, if $A\in
H_2(X_i;\mathbb{Z})$ is represented by a $J_i$-holomorphic sphere
for some almost complex structure $J_i$ preserving $TF_i$, then we
necessarily have $\langle\kappa_{X_i}+PD[F_i],A\rangle\geq 0$.

As such, if the symplectic sum is smoothly nontrivial and relatively
minimal, and if the genus of the surfaces involved is positive, then
(\ref{minusone}) cannot hold, and so $Z$ cannot admit any
nonvanishing genus zero Gromov--Witten invariants in any nontrivial
homology classes.  In \cite{U} this was used to conclude that $Z$ is
minimal; in our present case, we simply note that any rational or
ruled surface does admit a nonvanishing genus-zero Gromov--Witten
invariant in a nontrivial homology class (the proper transform of
the hyperplane class for blowups of $\mathbb{C}P^2$, or the class of
the fiber of a ruling for ruled surfaces; see, \emph{e.g.},
\cite{Mc1}). Thus rational and ruled surfaces (\emph{i.e.},
symplectic $4$-manifolds of Kodaira dimension $-\infty$) cannot
arise as smoothly nontrivial, relatively minimal symplectic sums
along surfaces of positive genus, proving part (a) of Theorem
\ref{main}.

\section{$4$-manifolds of Kodaira dimension zero}

The proof of part (b) of Theorem \ref{main} depends on the following
result.

\begin{theorem} \label{minusf} Let $(X_1,\omega_1),(X_2,\omega_2)$ be
symplectic $4$-manifolds, $F_i\subset X_i$ ($i=1,2$) embedded
symplectic submanifolds of equal positive genus, equal area, and
opposite self-intersection, and $(Z=X_1\#_{F_1=F_2} X_2,\omega)$ the
symplectic sum of the $X_i$ along the $F_i$.  Assume that the
symplectic sum is smoothly nontrivial and that $(Z,\omega)$ has
Kodaira dimension zero and is minimal. Then the $F_i$ are both tori,
and are Poincar\'e dual to $-\kappa_{X_i}$.\end{theorem}

\begin{proof} According to Theorem 2.1 of \cite{IP}, there is a
symplectic $6$-manifold $(\mathcal{Z},\Omega)$ equipped with a
projection $\pi\co \mathcal{Z}\to D^2$ with the property that, for
$\lambda\neq 0$, $(\pi^{-1}(\lambda),\Omega|_{\pi^{-1}(\lambda)})$
is isotopic to the symplectic sum $Z$, while
$(\pi^{-1}(0),\Omega|_{\pi^{-1}(0)})$ is the singular symplectic
manifold obtained by directly gluing the $X_i$ along the $F_i$, so
that the two pieces intersect in a copy $F$ of the $F_i$. As seen in
the proof of Lemma 2.2 of \cite{IP}, one has (for $\lambda\neq 0$)
$\kappa_{\mathcal{Z}}|_{\pi^{-1}(\lambda)}=\kappa_{\pi^{-1}(\lambda)}$,
while for $i=1,2$ $c_1(T\mathcal{Z}|_{X_i})=c_1(X_i)-PD_{X_i}[F_i]$.

We claim now that,  where $i_{\lambda}\co \pi^{-1}(\lambda)\to \mathcal{Z}$,
$i_j\co X_j\to \mathcal{Z}$ $(j=1,2$) are the inclusions,
$i_{\lambda\ast}[\pi^{-1}(\lambda)]$ is equal to
$i_{1\ast}[X_1]+i_{2\ast}[X_2]$ in $H_4(\mathcal{Z};\mathbb{Z})$.
This essentially follows from the description of $\mathcal{Z}$ in Section 2 of \cite{IP}.  Namely, as seen there,
 $\pi^{-1}(\lambda)$ for
$0<\lambda\ll 1$ may (up to isotopy) be obtained from
$\pi^{-1}(0)=X_1\cup_{F}X_2$ as follows: a neighborhood of
$F\subset\mathcal{Z}$ is diffeomorphic to a neighborhood $U$ of the
zero section in $L\oplus L^*\to F$ where $L\to F$ is a complex line
bundle of degree $[F_1]^2$ and $L^*$ is its dual; putting a Hermitian metric $|\!\cdot\!|$ on $L$ (which then induces one on $L^*$)
and arranging that $U=\{(\alpha,v,x)\in L\oplus L^*||\alpha|,|v|\leq 1\}$, up to smooth
isotopy $\pi^{-1}(\lambda)$ coincides with $\pi^{-1}(0)$ outside $U$
and its intersection with $U$ is given by $\{(\alpha,v,x)\in L\oplus
L^*|\langle \alpha,v\rangle=\chi(|\alpha|+|v|)\lambda\}$ where $\chi\co [0,1]\to [0,1]$ is
a bump function supported in $[0,2\lambda^{1/4}]$ and equal to $1$ on $[0,\lambda^{1/4}]$.
Observe that $X_1$ is homologous to the chain obtained by replacing $X_1\cap U$ by (in hopefully self-explanatory notation)  \[
\{(\bar{v},v,x)| |v|\leq \lambda^{1/2}\}+\{(\alpha,v,x)\in \pi^{-1}(\lambda)|\lambda^{1/2}\leq |v|\leq 1\}\] and likewise $X_2$ is homologous to the chain obtained by replacing $X_2\cap U$ by \[
\{(\alpha,\bar{\alpha},x)| |\alpha|\leq \lambda^{1/2}\}+\{(\alpha,v,x)\in \pi^{-1}(\lambda)|\lambda^{1/2}\leq |\alpha| \leq 1\}.\]  But adding these two chains together and passing to homology just recovers $i_{\lambda\ast}[Z_{\lambda}]$ (as the first terms in the two expressions above now cancel), verifying that indeed $i_{\lambda\ast}[\pi^{-1}(\lambda)]=i_{1\ast}[X_1]+i_{2\ast}[X_2]\in H_4(\mathcal{Z};\mathbb{Z})$.

As such,  we have
\begin{align} \langle\kappa_Z\cup[\omega],[Z]\rangle&=\langle
i_{\lambda}^{*}(\kappa_{\mathcal{Z}}\cup[\Omega]),[\pi^{-1}(\lambda)]\rangle=\langle\kappa_{\mathcal{Z}}\cup[\Omega],i_{\lambda\ast}[\pi^{-1}(\lambda)]\rangle\nonumber\\
&=\langle\kappa_{\mathcal{Z}}\cup[\Omega],i_{1\ast}[X_1]\rangle+\langle\kappa_{\mathcal{Z}}\cup[\Omega],i_{2\ast}[X_2]\rangle
\nonumber\\&=\langle\kappa_{\mathcal{Z}}|_{X_1}\cup
[\Omega]|_{X_1},[X_1]\rangle+\langle\kappa_{\mathcal{Z}}|_{X_2}\cup
[\Omega]|_{X_2},[X_2]\rangle\nonumber
\\&=\langle(\kappa_{X_1}+PD_{X_1}[F_1])\cup
[\omega_1],[X_1]\rangle+\langle(\kappa_{X_2}+PD_{X_2}[F_2])\cup[\omega_2],[X_2]\rangle
\nonumber\\&=\langle\kappa_{X_1}\cup[\omega_1],[X_1]\rangle+\langle\kappa_{X_2}\cup[\omega_2],[X_2]\rangle+\int_{F_1}\omega_1+\int_{F_2}\omega_2.\label{sumkappa}\end{align}
Note that nothing that we've done so far makes any use of the
assumption that $(Z,\omega)$ is minimal with Kodaira dimension zero.
However, when we do implement that assumption, the left hand side of
(\ref{sumkappa}) becomes zero.  Hence \[
\left(\langle\kappa_{X_1}\cup[\omega_1],[X_1]\rangle+\int_{F_1}\omega_1\right)+\left(\langle\kappa_{X_2}\cup[\omega_2],[X_2]\rangle+\int_{F_2}\omega_2\right)=0.\]
We claim now that each $\kappa_{X_i}$ is proportional to $PD[F_i]$
in $H^2(X_i;\mathbb{R})$.  Indeed, suppose that this were not the
case, say for $i=1$.  Then we could find an element $\beta\in
H^2(X_1;\mathbb{R})$ such that
$\langle\kappa_{X_1}\cup\beta,[X_1]\rangle=1$ but
$\langle\beta,[F_1]\rangle=0$.  Moreover, since small closed
perturbations of symplectic forms are still symplectic, for any
sufficiently small $\epsilon>0$ the cohomology class
$[\omega_1]+\epsilon\beta$ would admit a symplectic form, say
$\omega_{1}^{\epsilon}$, which (assuming $\epsilon$ is small
enough) would continue to make $F_1$ symplectic, and would induce the same canonical class $\kappa_{X_1}$ as does $\omega_1$ by virtue of being deformation equivalent to $\omega_1$.  We would then
have \[
\int_{F_1}\omega_{1}^{\epsilon}=\int_{F_1}\omega_1,\quad\langle\kappa_{X_1}\cup[\omega_{1}^{\epsilon}],[X_1]\rangle=\langle\kappa_{X_1}\cup[\omega_1],[X_1]\rangle+\epsilon.\]
But then we could apply the symplectic sum operation to the
symplectic manifolds $(X_1,\omega_{1}^{\epsilon})$,
$(X_2,\omega_{2})$ along the $F_i$ to obtain a new symplectic
manifold $(Z^{\epsilon},\omega^{\epsilon})$, diffeomorphic (indeed
deformation equivalent) to $(Z,\omega)$, and by (\ref{sumkappa})
we would have \[
\langle\kappa_{Z^{\epsilon}}\cup[\omega],[Z^{\epsilon}]\rangle=\epsilon>0.\]
So $(Z^{\epsilon},\omega^{\epsilon})$ is not a minimal symplectic
manifold of Kodaira dimension zero.  But by Theorem 2.6 of
\cite{Li}, the Kodaira dimension of a symplectic four-manifold is a
diffeomorphism invariant, so since $Z^{\epsilon}$ is diffeomorphic
to $Z$ (and since the minimality or nonminimality of $Z$ is also a
diffeomorphism invariant; see, \emph{e.g.}, Proposition 2.1 of
\cite{Li}) this is impossible. This contradiction shows that, for
some $\mu_1,\mu_2\in\mathbb{R}$, we have $PD(\kappa_{X_i})=\mu_i[F_i]\in
H_2(X;\mathbb{R})$ for $i=1,2$.

Evidently, since $\langle\kappa_Z\cup[\omega],[Z]\rangle=0$, we must
have \[
(2+\mu_1+\mu_2)\int_{F_1}\omega_1=\left(\langle\kappa_{X_1}\cup[\omega_1],[X_1]\rangle+\int_{F_1}\omega_1\right)+\left(\langle\kappa_{X_2}\cup[\omega_2],[X_2]\rangle+\int_{F_2}\omega_2\right)=0,\]
so since $\int_{F_1}\omega_1> 0$ we have \[ \mu_1+\mu_2=-2.\]  In
particular, at least one of the $\mu_i$, say $\mu_1$, is negative.  Then
$\langle\kappa_{X_1}\cup[\omega_1],[X_1]\rangle=\mu_1\int_{F_1}\omega_1<0$,
so by Theorem B of \cite{Liu} and results of \cite{T}
$(X_1,\omega_1)$ must be either a rational surface or an irrational
ruled surface.

Now if $X_1=\mathbb{C}P^2\#k\overline{\mathbb{C}P^2}$, then where
$H$ is the hyperplane class and $E_1,\ldots,E_k$ are the classes of
the exceptional spheres of the blowups, we have
$PD(\kappa_{X_1})=-3H+\sum_{i=1}^{k}E_i$; from this it's easy to see
from the adjunction formula that if $F_1$ is a symplectic surface of
positive genus representing $\frac{1}{\mu_1}PD(\kappa_X)$ we
necessarily have $-\frac{1}{\mu_1}\geq 1$, so that $\mu_1\geq -1$;
similarly if $X_1=S^2\times S^2$ then
$PD(\kappa_{X_1})=-2[S^2\times\{pt\}]-2[\{pt\}\times S^2]$, and
again our assumptions on $F_1$ force $\mu_1\geq -1$ ($\mu_1=-2$ is ruled
out both by the assumption that $F_1$ has positive genus and by the
assumption that the sum is smoothly nontrivial, so that $F_1$ is not
a section of an $S^2$-bundle). Thus if $X_1$ is rational then
$\mu_1\geq -1$.

Suppose now that $X_1$ is  an irrational ruled surface; say
$X_1=(S^2\times \Sigma_h)\#k\overline{\mathbb{C}P^2}$ $(h\geq
1,k\geq 0$). Then where $\sigma$ is the homology class of the proper
transform of a section of $S^2\times \Sigma_h\to \Sigma_h$ and $f$
is the class of the proper transform of a generic fiber, we have
$PD(\kappa_{X_1})=-2\sigma+(2h-2)f+\sum_{i=1}^{k}e_i$.  If $k>0$,
then in order for $[F_1]$ to be an integral class we necessarily
have $\mu_1\geq -1$.  If $k=0$, we have either $\mu_1\geq -1$ or
$\mu_1=-2$, but if $\mu_1=-2$ then by Proposition 3.3 of \cite{U} there
is a ruling on $X_1$ of which $F_1$ is a section, which is forbidden
by the assumption that the sum is smoothly nontrivial.  Finally, if
$X_1$ is a nontrivial $S^2$-bundle $S^2\tilde{\times}\Sigma_h$ over
a surface of genus $h>0$, then where $s^+$, $s^-$ are sections of
square $+1$ and $-1$ respectively, we have
$PD(\kappa_X)=(2h-3)[s^+]-(2h-1)[s^-]$, so that since
$[F_1]=\frac{1}{\mu_1}PD(\kappa_X)\in H_2(X_1;\mathbb{Z})$, we get
$|\frac{1}{\mu_{1}}|\geq 1$, so $\mu_1\geq -1$.

Summing up, we've shown that since $\mu_1<0$, $X_1$ is rational or
ruled, and that (because $X_1$ is rational or ruled) by the hypothesis of the theorem we have $\mu_1\geq
-1$.  But then since $\mu_1+\mu_2 = -2$, we must have $\mu_2<0$, which
then forces $X_2$ to be rational or ruled, and so by the same
arguments $\mu_2\geq -1$.  So since $\mu_1+\mu_2=-2$ and $\mu_1,\mu_2\geq -1$,
we in fact have $\mu_1=\mu_2=-1$, \emph{i.e.}, $[F_i]=-PD(\kappa_{X_i})$
for $i=1,2$.

To see that the $F_i$ are tori, simply note that the adjunction
formula gives \[
2g(F_i)-2=[F_i]^2+\langle\kappa_{X_i},[F_i]\rangle=\kappa_{X_i}^2-\kappa_{X_i}^{2}=0.\]

\end{proof}

Theorem \ref{minusf} and its proof show that, in order for a
smoothly nontrivial symplectic sum $X_1\#_{F_1=F_2} X_2$ along a
surface of positive genus to be minimal and of Kodaira dimension
zero, the $X_i$ necessarily  are rational
or ruled, and the $F_i$ are tori Poincar\'e dual to $-\kappa_{X_i}$.
The main theorem of \cite{U} shows that the sum will automatically
be minimal provided that the $X_i$ and $F_i$ satisfy the hypotheses
of Theorem \ref{main}. The problem of determining the diffeomorphism
types of such sums now naturally splits into 3 cases, which we
address in turn.

\subsection{Case 1: $X_1$ and $X_2$ are rational} Note that, since
the $F_i$ are tori, we have $\kappa_{X_1\#_{F_1=F_2}
X_2}^{2}=3\sigma(X_1\#_{F_1=F_2} X_2)+2\chi(X_1\#_{F_1=F_2}
X_2)=\kappa_{X_1}^{2}+\kappa_{X_2}^{2}$.  Now the rational surfaces
have $\kappa_{\mathbb{C}P^2\#k\overline{\mathbb{C}P^2}}^{2}=9-k$ and
$\kappa_{S^2\times S^2}^{2}=8$, in light of which, in order for
$\kappa_{X_1\#_{F_1=F_2} X_2}^{2}$ to be zero, the unordered pair
$\{X_1,X_2\}$ must be either
$\{\mathbb{C}P^2\#k\overline{\mathbb{C}P^2},\mathbb{C}P^2\#(18-k)\overline{\mathbb{C}P^2}\}$
($0\leq k\leq 9$) or $\{S^2\times
S^2,\mathbb{C}P^2\#17\overline{\mathbb{C}P^2}\}$. Note that since
$[F_i]=-PD(\kappa_{X_i})$, $F_i$ has intersection number $1$ with
each embedded symplectic $(-1)$-sphere in $X_i$.  Choosing an almost
complex structure $J_i$  on $X_i$ generic among those making $F_i$
pseudoholomorphic, each homology class in $X_i$ that is represented
by an embedded symplectic $(-1)$-sphere will be represented by a
$J_i$-holomorphic $(-1)$-sphere, which will then intersect $F_i$
transversely and once.  Blowing down such a $(-1)$-sphere will then
result in a symplectic manifold $X'_i$ together with an embedded
torus $F'_i$, with $[F'_i]^2=[F_i]^2+1$.

Now according to Lemma 5.1 of \cite{G} (respectively, Proposition
1.6 of \cite{MSy}), if $(M_1,\omega_1),(M_2,\omega_2)$ are
symplectic $4$-manifolds containing surfaces $\Sigma_1,\Sigma_2$ of
equal genus and area such that $[\Sigma_1]^2+[\Sigma_2]^2=1$, and if
$(\tilde{M_i},\tilde{\Sigma}_i)$ denotes the pair consisting of a
symplectic $4$-manifold and symplectic surface that results from
blowing up $M_i$ at a point of $\Sigma_i$, then
$\tilde{M_1}\#_{\tilde{\Sigma}_1=\Sigma_2}M_2$ is diffeomorphic
(respectively, symplectic deformation equivalent) to
$M_{1}\#_{\Sigma_1=\tilde{\Sigma}_2}\tilde{M_2}$.

Applying this to the cases under consideration enables us to replace
the pair
$\{\mathbb{C}P^2\#k\overline{\mathbb{C}P^2},\mathbb{C}P^2\#(18-k)\overline{\mathbb{C}P^2}\}$
$(0\leq k\leq 9$) by
$\{\mathbb{C}P^2\#(k+1)\overline{\mathbb{C}P^2},\mathbb{C}P^2\#(18-k-1)\overline{\mathbb{C}P^2}\}$
(since if $X_2=\mathbb{C}P^2\#(18-k)$ where $k\leq 9$, $X_2$
contains $(-1)$-spheres meeting $F_2$ transversely once, any one of
which may be blown down to give a torus of square one larger in
$\mathbb{C}P^2\#(18-k-1)\overline{\mathbb{C}P^2}$, and we may then
apply Gompf's and McDuff-Symington's results mentioned above).
Repeatedly ``trading blowups'' in this fashion enables us to reduce
to the case that $k=18-k=9$. Similarly, if $X_1=S^2\times S^2$, then
$X_2=\mathbb{C}P^2\#17\overline{\mathbb{C}P^2}$, so ``trading
blowups'' once reduces us to the case that $X_1=(S^2\times
S^2)\#\overline{\mathbb{C}P^2}=\mathbb{C}P^2\#2\overline{\mathbb{C}P^2}$
and $X_2=\mathbb{C}P^2\#16\overline{\mathbb{C}P^2}$, and then doing
so seven more times reduces to the case that
$X_1=X_2=\mathbb{C}P^2\#9\overline{\mathbb{C}P^2}$.

So we assume for the rest of this subsection that
$X_1=X_2=\mathbb{C}P^2\#9\overline{\mathbb{C}P^2}$.  The $X_i$ are
then both diffeomorphic to the total space of a rational elliptic
surface $E(1)$, and the $F_i$ are homologous to a fiber of $E(1)$.
In fact, we can make a stronger statement about the $F_i$; we begin
with the following lemma.

\begin{lemma} \label{makehol} Let $B\subset \mathbb{C}^2$ be a closed ball around the origin, and let
$\phi\co B\to\mathbb{C}^2$ be an orientation preserving
diffeomorphism onto a closed subset of $\mathbb{C}^2$ which maps
$\{(w,z)\in B|z=0\}$ orientation-preservingly to $\mathbb{C}\times
\{0\}$.  Then there is a diffeomorphism $\phi'\co B\to \phi(B)$
which is isotopic to $\phi$ by an isotopy which is the identity on a
neighborhood of $\partial B$, such that $\phi'$ maps $\{(w,z)\in
B|z=0\}$ to $\mathbb{C}\times\{0\}$ and is holomorphic on a
neighborhood of the origin.
\end{lemma}

\begin{proof}  With respect to the splitting
$\mathbb{R}^4=\mathbb{C}^2=\mathbb{R}^2\times\mathbb{R}^2$, the
Jacobian $J$ of $\phi$ at the origin has the $2\times 2$ block form
$\left(\begin{array}{ll}A&M\\0&C\end{array}\right)$. The
orientation-preserving conditions ensure that $A$ and $J$ both have
positive determinant, in light of which $C$ also has positive
determinant.  Let $\chi\co B\to [0,1]$ be a smooth function which is
equal to $1$ on a neighborhood of the origin and to $0$ on a
neighborhood of $\partial B$. Let $V_1$ be the vector field on
$B\subset \mathbb{R}^4$ whose value at $v\in B$ is $\chi(v)
\left(\begin{array}{cc}-\log A&0\\0&-\log C\end{array}\right)v$; by
composing $\phi$ with the time-$1$ map of $V_1$ we reduce to the
case where the Jacobian at the origin has form
$\left(\begin{array}{cc}I&A^{-1}M\\0&I\end{array}\right)$. Then
where
$V_2(v)=\chi(v)\left(\begin{array}{cc}0&-A^{-1}M\\0&0\end{array}\right)v$
composing with the time-$1$ map of $V_2$ reduces us to the case that
the Jacobian of $\phi$ at the origin is the identity.  (Note that
both $V_1$ and $V_2$ are tangent to $\mathbb{C}\times\{0\}$, so the
condition that this set be preserved is not disturbed).  But in this
case there is a constant $C$ (depending on the second derivatives of
$\phi$) such we have $|\phi(w,z)-(w,z)|\leq C(|w|^2+|z|^2)$ and
$|(D\phi)_{(w,z)}-I|\leq C(|w|^2+|z|^2)^{1/2}$.  Let $\beta\co
[0,\infty)\to [0,1]$ be a smooth, monotone function such that
$\beta(r)=0$ for $r\leq 1/3$, $\beta(r)=1$ for $r\geq 1$, and
$\beta'(r)\leq 2$.  Then if $\delta>0$ defining \[
\phi'(w,z)=(w,z)+\beta((|w|^2+|z|^2)^{1/2}/\delta)(\phi(w,z)-(w,z))\]
yields a smooth map which coincides with $\phi$ outside the ball of
radius $\delta$ around the origin, is holomorphic in the ball of
radius $\delta/3$ around the origin, and differs from $\phi$ in
$C^1$ norm by at most $(3+\delta)C\delta$ and therefore is a
diffeomorphism as long as $\delta$ is taken small enough.  Also,
since $\phi$ preserves $\{z=0\}$, $\phi'$ evidently does as well,
and $\phi'$ is isotopic to $\phi$ by the isotopy
$\phi+t(\phi'-\phi)$.

\end{proof}

\begin{lemma}\label{isotdiff}
Let $(M,\omega)$ be a symplectic manifold obtained by blowing up
$(\mathbb{C}P^2,\omega_{FS})$ at nine distinct points, and let
$F\subset M$ be an embedded symplectic submanifold Poincar\'e dual to the
anticanonical class $-\kappa_M$.  Then there exists an elliptic
fibration $\pi\co E(1)\to S^2$ with no multiple fibers and a
diffeomorphism $\Phi\co M\to E(1)$ with the property that $\Phi(F)$
is a fiber of $\pi$.
\end{lemma}
\begin{proof}
Let $\tilde{J}_0$ be an almost complex structure which makes $F$
pseudoholomorphic.  Provided that $\tilde{J}_0$ is chosen
generically from among such almost complex structures, there are
unique, embedded $\tilde{J}_0$-holomorphic representatives
$E_1,\ldots,E_9$, respectively, of the homology classes of the nine
exceptional divisors of the blowup $M\to\mathbb{C}P^2$.  Each of the
$E_i$ has a tubular neighborhood $N_i$ modeled (symplectically) on a
neighborhood of the zero section of the holomorphic line bundle
$\mathcal{O}(-1)$ over $\mathbb{C}P^1$; we may shrink the $N_i$ if
necessary to make them pairwise disjoint.  Now, for each $i$, the
$\tilde{J}_0$-holomorphic curves $F$ and $E_i$ have homological
intersection number $1$, so they meet transversely, positively, and
in just one point.  Hence as in Lemma 2.3 of \cite{G}, we may
isotope $F$ rel $M\setminus\cup N_i$ by symplectomorphisms to a
surface $F'$ such that, for some smaller neighborhoods $N'_i$ of
$E_i$, $F'\cap N'_i$ coincides with a fiber of the bundle projection
$\mathcal{O}(-1)\to E_i$.  Further, during this isotopy, we can
change the almost complex structure $\tilde{J}_0$ to an almost
complex structure $\tilde{J}$ which agrees with $\tilde{J}_0$
outside $N_i$, makes $F'$ holomorphic, and coincides on $N'_i$ with
the standard integrable complex structure on the tautological line
bundle $\mathcal{O}(-1)$.

We now symplectically blow down the $E_i$ (see Section 7.1 of
\cite{McDS} for a detailed description of this process); doing so
amounts to symplectically identifying an annular neighborhood
$L_i=D(\delta_{1,i})\setminus D(\delta_{2,i})\subset N'_i\subset
\mathcal{O}(-1)$ of the zero section $E_i$ with a spherical shell
$B(r_{1,i})\setminus B(r_{2,i})$ in $\mathbb{C}^2$, and filling in
this shell by a standard ball (centered at a point $p_i$ with radius
$r_{2,i}$) to replace the $\delta_{2,i}$-neighborhood of $E_i$. This
process is compatible with the complex blowdown in the sense that
the annulus fibers of $L_i$ are taken to annuli in the complex lines
that they correspond to under the identification of
$\mathcal{O}(-1)$ with the tautological line bundle of
$\mathbb{C}^2$, so in the blowdown these annuli may be filled in to
form discs by adding in the disc of radius $r_{2,i}$ in the
corresponding complex line.  So since $F'$ meets each $L_i$ in one
of these annulus fibers, filling these annuli in as above gives rise
to a compact symplectic surface $S$ whose intersection with each
$B(r_{1,i})$ coincides with a complex line.

As a result of all this, blowing down the $E_i$ results in
$\mathbb{C}P^2$ equipped with an almost complex structure $J$ which
is integrable near each of the points $p_i$ (and coincides with
$\tilde{J}$ outside the blow-up neighborhoods), a $J$-holomorphic
curve $S\subset\mathbb{C}P^2$ which coincides with $F'$ outside the blow-up
neighborhoods, and a symplectic form $\omega'$ (induced by $\omega$
and the blowdown procedure) which is compatible with $J$. $\omega'$
is easily seen to be cohomologous to the Fubini-Study form (for
instance, the proper transform of the hyperplane class under the
original blowups that were done to obtain $M$ from $\mathbb{C}P^2$
has a nonvanishing Gromov--Witten invariant, and so may be
represented by a symplectic surface $\tilde{H}$ which we can arrange
to miss the $L_i$ and so gives rise to a symplectic surface in
$(\mathbb{C}P^2,\omega')$ with the same area as $\tilde{H}$ has in
$(M,\omega)$).  Hence by a result of \cite{Gr} $\omega'$ is
symplectomorphic to the Fubini-Study form.  Also, for a similar
reason as above, $S\subset\mathbb{C}P^2$ has degree 3.  So by
Theorem 3 of \cite{Sik}, $S$ is symplectically isotopic to some
complex cubic curve $C$.  The isotopy extension theorem then ensures
that there is an ambient isotopy
$\phi_t:\mathbb{C}P^2\to\mathbb{C}P^2$ such that $\phi_1(S)=C$.  Our
intention now is to modify $\phi_1$ to some other diffeomorphism
$\phi_2$  so that (i) $\phi_2(S)=C$, (ii) $\mathbb{C}P^2$ carries an
elliptic pencil with fiber $C$ and base locus
$\{\phi_2(p_1),\ldots,\phi_2(p_9)\}$, and (iii) $\phi_2$ lifts to a
diffeomorphism of the appropriate blowups which takes $F'$ to the
proper transform of $C$.

To achieve this, first note that since the sheaf of sections of the
holomorphic normal bundle to $C$ which vanish at the eight points
$\phi_1(p_1),\cdots,\phi_1(p_8)$ has degree $1$, it admits global
nonvanishing sections by Riemann-Roch; perturbing $\phi_1$ by
composing with a small diffeomorphism which preserves $C$ we can
arrange that $\phi_1(p_1),\cdots,\phi_1(p_8)$ are generic in the
sense that none of these sections vanishes to order $2$.  Such a
section results in another smooth cubic $C'$ which meets $C$
transversely at the eight points $\phi_1(p_1),\ldots,\phi_1(p_8)$,
and also at some other point $q_9$.  Now let $v$ be a vector field
on $C$ which vanishes on neighborhoods of the $\phi_1(p_i)$ for
$1\leq i\leq 8$ and whose time-$1$ flow maps $\phi_1(p_9)$ to $q_9$.
Using a partition of unity subordinate to a set of local
trivializations of the normal bundle to $C$ and a bump function
supported on a tubular neighborhood of $C$, extend $v$ to a vector
field $V$ on $\mathbb{C}P^2$ whose restriction to $C$ is $v$; let
$\psi$ be the time-1 flow of $V$, and let $\phi_0=\psi\circ\phi_1$.
Then $\phi_0\co\mathbb{C}P^2\to\mathbb{C}P^2$ is a diffeomorphism
which maps $S$ to $C$ and $p_1,\ldots,p_9\in S$ to the $9$
intersection points between $C$ and another smooth cubic $C'$.

Now $\phi_0$ will not lift to a diffeomorphism on blowups, because
it's not holomorphic near the points being blown up and so doesn't
map complex lines to complex lines.  However, around each $p_i$
there are complex coordinates $(w,z)$ in which $S$ is given by
$\{z=0\}$ and likewise near each $\phi_0(p_i)$ there are complex
coordinates in which $C$ is given by $\{z=0\}$.  Hence, in terms of
these local holomorphic coordinates, $\phi_0$ is given in these
neighborhoods by a diffeomorphism which satisfies the hypothesis of
Lemma \ref{makehol}.  So $\phi_0$ may be modified by an isotopy
supported in the union of these neighborhoods to an
orientation-preserving diffeomorphism
$\phi_2\co\mathbb{C}P^2\to\mathbb{C}P^2$ which maps $S$ to $C$ and
is holomorphic on (smaller) neighborhoods of $p_1,\ldots,p_9$.
Consequently, if $Y$ is the (complex) blowup of $(\mathbb{C}P^2,J)$
at the 9 points $p_1,\ldots,p_9$ (which makes sense because
$J$ is integrable near the $p_i$), and if $E(1)$ is the complex
blowup of $(\mathbb{C}P^2,J_{std})$ at
$\phi_2(p_1),\ldots,\phi_2(p_9)$, $\phi_2$ lifts to a diffeomorphism
$Y\to E(1)$ taking the proper transform of $S$ to the proper
transform of $C$.  If $f$ and $f'$ are homogeneous cubic polynomials
with vanishing loci $C$ and $C'$ respectively, the vanishing loci
$C_{[\lambda:\mu]}$ of $\lambda f+\mu f'$
$([\lambda:\mu]\in\mathbb{C}P^1$) provide an elliptic pencil on
$\mathbb{C}P^2$ with base locus $\{\phi(p_1),\ldots,\phi(p_9)\}$.
Blowing up the base locus to form $E(1)$ thus gives an elliptic
fibration with the proper transform of $C$ as a fiber.

To compare $Y$ to $X$, note that in order to put a symplectic form
on the complex blowup $Y$ in such a way that the map $\phi_2$ lifts
to $Y$ we need to cut out balls $B'_i$ around $p_i$ that are smaller
than the balls $B(r_{2,i})$ that were created by the blowdown
$X\to\mathbb{C}P^2$. $X$ and $Y$ hence can't be symplectomorphic.
However, there is an obvious diffeomorphism between the blowups
corresponding to balls of different size which simply changes the
radius of the disc bundle involved; in particular this
diffeomorphism has a restriction to the neighborhood of the
exceptional sphere which preserves the fibers of the normal bundle
$\mathcal{O}(-1)$. Recalling that $F'$ coincided with a fiber of the
normal bundle on the neighborhood $N'_i$ it then follows that the
natural diffeomorphism $X\to Y$ takes $F'$ to the proper transform
of $S$ in $Y$.

We hence have a diffeomorphism $X\to E(1)$ which takes $F'$ to the
proper transform of $C$; precomposing this with a symplectomorphism
of $X$ which isotopes $F$ to $F'$ gives us the promised map $\Phi\co
X\to E(1)$ taking $F$ to a fiber of an elliptic fibration.
\end{proof}

\begin{cor} If the symplectic sum $Z$ of two rational surfaces along a positive genus surface is a
minimal symplectic four-manifold of Kodaira dimension zero, then $Z$
is diffeomorphic to the $K3$ surface.
\end{cor}

\begin{proof}  Denote the two summands by $X_i$ ($i=1,2$) and the surfaces in
question by $F_i$.  We've shown that we may assume that
$X_i=\mathbb{C}P^2\#9\overline{\mathbb{C}P^2}$ and that the $F_i$
are tori Poincar\'e dual to $-PD(\kappa_{X_i})$.  Hence by Lemma
\ref{isotdiff}, we have $X_i\setminus \nu F_i\cong E(1)\setminus \nu
F$ where $E(1)$ is the total space of a rational elliptic fibration
having fiber $F$ (and where ``$\cong$'' denotes diffeomorphism).  Hence \[ Z\cong (E(1)\setminus \nu
F)\cup_{\Phi}(E(1)\setminus \nu F)\] for some orientation-reversing
diffeomorphism $\Phi$ of the boundary $\partial(E(1)\setminus \nu
F)=T^3$.  But according to Proposition 1 of Appendice 4 of
\cite{GM}, every orientation preserving diffeomorphism of
$\partial(E(1)\setminus \nu F)$ extends to $E(1)\setminus \nu F$, so
the diffeomorphism type of $(E(1)\setminus \nu
F)\cup_{\Phi}(E(1)\setminus \nu F)$ is independent of $\Phi$.  So
since one choice of $\Phi$ (namely the one corresponding to taking
the standard fiber sum of $E(1)$ with itself) gives rise to the $K3$
surface, $Z$ is evidently diffeomorphic to the $K3$ surface
independently of $\Phi$.
\end{proof}

\subsection{Case 2: $X_1$ is rational and $X_2$ is irrational and ruled}

Assume that $X_2$ is  a ruled surface over a curve $C$ of positive
genus; we'll show shortly that $C$ is a torus in the cases of
interest. Then $X_2$ is symplectomorphic either to the nontrivial
$S^2$-bundle over $C$, which we denote $S^2\tilde{\times}C$, or else
$X_2=(S^2\times C)\#k\overline{\mathbb{C}P^2}$ for some $k\geq 0$.
Suppose that $X_1$ is rational and $F_i\subset X_i$ ($i=1,2$) are
embedded symplectic submanifolds with the property that the smoothly
nontrivial symplectic sum $X_1\#_{F_1=F_2}X_2$ is minimal and of
Kodaira dimension zero. Then Theorem \ref{minusf} shows that the
$F_i$ are tori Poincar\'e dual to $-\kappa_{F_i}$. In particular, if
$X_2$ is nonminimal and $J_2$ is an almost complex structure
preserving $TF_2$, then each member of a maximal disjoint collection
of embedded $J_2$-holomorphic $(-1)$-spheres meets $F_2$
transversely and once; hence the results of \cite{G} and \cite{MSy}
alluded to in the previous subsection show that, up to deformation
equivalence, the symplectic sum is left unchanged if we
simultaneously blow down each member of this maximal collection and
blow up $X_1$ at a corresponding number of points on $F_1$.  This
reduces us to the case that $X_2$ is minimal, and so is either
$S^2\times C$ or $S^2\tilde{\times}C$.

\begin{lemma}\label{annbundle} Let $\pi\co E\to C$ be an $S^2$-bundle over a positive-genus surface $C$ with symplectic form $\omega\in \Omega^2(E)$,
and let $\Sigma\subset E$ be an embedded, connected, symplectic
representative of $-PD(\kappa_E)\in H_2(E;\mathbb{Z})$, with tubular
neighborhood $\nu\Sigma$.  Then $C$ is a torus and there is another bundle map $\pi'\co E\to C$ whose fibers are symplectic spheres such that $\pi'|_{E\setminus \nu\Sigma}$ defines a fiber bundle with fibers diffeomorphic to the
annulus $S^1\times I.$\end{lemma}

\begin{proof}
First, note that if $J$ is an almost complex structure compatible
with $\omega$ with respect to which $\Sigma$ is pseudoholomorphic,
and if $\mathcal{M}$ is the moduli space of unparametrized
pseudoholomorphic spheres representing the class of the fiber of
$\pi$, then results of \cite{Mc1} show that the map $\pi'\co
E\to\mathcal{M}$ which takes  $e\in E$ to the point of $\mathcal{M}$
representing the unique $J$-holomorphic representative of the fiber
class on which $e$ lies is an $S^2$-bundle with fibers homologous to
the fibers of $\pi$; fundamental group considerations then imply
that $\mathcal{M}$ has the same genus as $C$. By construction the fibers of $\pi'$ are $J$-holomorphic and hence symplectic.

We now claim that $\Sigma$ is isotopic to some surface $\Sigma'\subset E$ such that there is an almost complex structure $J'$ which makes $\Sigma'$ pseudoholomorphic and with respect to which $\pi'\co E\to\mathcal{M}$ is a pseudoholomorphic map (with respect to some almost complex structure on $\mathcal{M}$).   Indeed, as in the proof of Lemma 3.2 of \cite{LU}, using the parametrized Riemann mapping theorem we can find complex coordinates $(z,w)$ on a suitable open set $U_i$ of $E$ centered around any critical point $p_i$ of $\pi'|_{\Sigma}$ in terms of which the projection $\pi'\co E\to\mathcal{M}$ is given by $(z,w)\mapsto w$ and $\partial_{\bar{z}}$ lies in the $J$-antiholomorphic tangent space $T^{0,1}_{J}$.  The intersection of $\Sigma$ with this neighborhood will then be given by $\Sigma\cap U_i=\{w=g_i(z)\}$ where $g_i(z)=c_iz^{k_i}+O(|z|^{k_i+1})$ and $c_i\neq 0$; that $\pi'|_{\Sigma}$ has a critical point at $(0,0)$ amounts to the statement that $k_i>1$.  Note that since $\Sigma$ has intersection number $2$ with the fibers of $\pi'$ (as $[\Sigma]=-\kappa_E$ and the fibers are square-zero spheres), we in fact have $k_i=2$, and moreover there can only be one critical point of $\pi'|_{\Sigma}$ in any given fiber of $\pi'$; accordingly we can and do choose the $U_i$ so that the $\pi'(U_i)$ are disjoint as $i$ varies.  We can then use a cutoff function supported in $U_i$ and equal to $1$ on some smaller neighborhood $U'_i$ of the critical point $p_i$ to isotope $\Sigma$ rel $U_i$ to some new surface $\Sigma'$ whose intersection with $U'_i$ is given by $\Sigma'\cap U'_i=\{w=c_iz^{k_i}\}$; further, using the same cutoff, we can isotope $J$ rel $U_i$ to a new almost complex structure $J'$ which coincides with the standard integrable complex structure with holomorphic coordinates $(z,w)$ on $U'_i$, and with respect to which both $\Sigma'$ and the fibers of $\pi'$ are $J'$-holomorphic.  Repeating this near every critical point of $\pi'|_{\Sigma}$, $\Sigma'$ and the fibers of $\pi'$ are now $J'$-holomorphic and $\pi'\co E\to \mathcal{M}$ restricts to a pseudoholomorphic map on a neighborhood of $Crit(\pi'|_{\Sigma'})$ with respect to $J'$ on $E$ and a suitable complex structure on $\cup_i\pi'(U'_i)\subset \mathcal{M}$.  But then extending this complex structure to all of $\mathcal{M}$, since $\pi'|_{\Sigma\setminus \cup_iU'_i}$ is an unbranched cover a simple patching argument may be used to further modify $J'$ so that it continues to make $\Sigma'$ pseudoholomorphic and now also makes the whole projection $\pi'\co E\to\mathcal{M}$ pseudoholomorphic. This proves the claim at the start of this paragraph.

But now $\pi'|_{\Sigma'}\co \Sigma'\to \mathcal{M}$ is a holomorphic map from a torus to $\mathcal{M}$ with degree $2$; we know that $\mathcal{M}$ has the same positive genus as $C$, so it follows that that genus is one.  Hurwitz's formula then implies that $\pi'|_{\Sigma'}$ has no critical points; since the critical points of $\pi'|_{\Sigma'}$ were constructed to be just the same as those of $\pi'|_{\Sigma}$ it then follows that $\pi'|_{\Sigma}$ has no critical points.
Thus $\Sigma$ meets every fiber of $\pi'$ transversely, and
hence exactly twice by the positivity of intersections of $J$-holomorphic curves.  Hence \[
\pi'|_{E\setminus\nu\Sigma}\co E\setminus \nu\Sigma\to T^2\] is a
fibration with fiber given by the complement of two discs in $S^2$,
\emph{i.e.}, by $S^1\times I$.
\end{proof}

Consequently, in all cases of interest, we have $X_2=S^2\times T^2$
or $X_2=S^2\tilde{\times}T^2$; these both have $c_{1}^{2}=0$, and so
if $X_1\#_{F_1=F_2}X_2$ is to be minimal of Kodaira dimension zero
then $c_{1}^{2}(X_1)=0$.  Since the only rational surface with
$c_{1}^{2}=0$ is $\mathbb{C}P^2\#9\overline{\mathbb{C}P^2}$,
evidently $X_1=\mathbb{C}P^2\#9\overline{\mathbb{C}P^2}$.

The above lemma makes the diffeomorphism classification of annulus
bundles over $T^2$ relevant to us; specifically we are interested in
those annulus bundles with orientable total space and having just
one boundary component.  Identify $S^1\times I$ with
$A=\overline{D(2)}\setminus D(1/2)\subset\mathbb{C}$. Any annulus bundle over
$T^2$ is isomorphic to one of the form $M(f,g;\{h_t\})$ where
$f,g\in \pi_0(Diff(A))$ commute, $\{h_t\}_{t\in S^1}\in
\pi_1(Diff(A))$, \[ M(f,g;\{1\})=\frac{R^2\times A}{(x+1,y,z)\sim
(x,y,f(z)),(x,y+1,z)\sim(x,y,g(z))}\] and $M(f,g;\{h_t\})$ is
obtained from $M(f,g;\{1\})$ by removing a trivial neighborhood
$D^2\times A$ from $M(f,g;\{1\})$ and gluing it back by
the map \begin{align*} \partial D^2\times A &\to \partial D^2\times
A\\ (t,z)&\to (t,h_t(z)).\end{align*} Since changing the choice of
basis $\{u,v\}$ of $H_1(T^2;\mathbb{Z})$ to, respectively,
$\{u+v,v\}$ or $\{v,u\}$ corresponds to replacing $(f,g)$ by
$(f\circ g,g)$ or $(g,f)$, we can assume that $f$ maps each
respective boundary component of $A$ to itself (if $f$ doesn't
initially, then either $g$ or $f\circ g$ does).  Now
$\pi_0(Diff(A))=\mathbb{Z}_2\oplus \mathbb{Z}_2$, with generators
given by $z\mapsto z^{-1}$ and $z\mapsto \bar{z}$.

We're interested in orientable annulus bundles over $T^2$ having
just one boundary component.  The orientability condition restricts
us to the case that the monodromies $f$ and $g$ preserve the
orientation of $A=\overline{D(2)}\setminus D(1/2)$. $f$ is assumed to map each
boundary component to itself, so this forces $f$ to be isotopic to
the identity.  But then in order for the bundle to have just one
boundary component $g$ must swap the boundary components of $A$, forcing
$g$ to be isotopic to $z\mapsto z^{-1}$.

Now as explained after the statement of Theorem 2.3 in \cite{Yag},
it follows from a theorem of Smale that the identity component of
$Diff(A)$ retracts to $S^1$ (and indeed the map $Diff_0(A)\to
Diff_0(S^1)$ given by restriction to one boundary component is a
homotopy equivalence), so $\pi_1(Diff(A))$ is generated by the loop
of diffeomorphisms $r_t\co A\to A$ where for $t\in S^1$ $r_t$ is
given by rotation through the angle $t$.  Thus, any orientable
annulus bundle over $T^2$ with one boundary component has form
\[ Y_n=M(I,z\mapsto z^{-1};\{r_{t}^{n}\})\] In fact, arguing exactly as in
Lemma 7 of Section 8 of \cite{S}, where $l,m\in \pi_1(\partial
D^2\times A)$ are, respectively, the generators of the images of the
inclusion-induced maps $\pi_1(A)\to \pi_1(\partial D^2\times A)$ and
$\pi_1(\partial D^2)\to \pi_1(\partial D^2\times A)$, one finds that
for any $n\in\mathbb{Z}$ there is a fiber preserving diffeomorphism
$Y_0\setminus D^2\times A\to Y_0\setminus D^2\times A$ whose
restriction to the boundary $\partial D^2\times A$ takes a
representative of $m$ to a representative of $m+2nl$ (an explicit formula for such a diffeomorphism may easily be found by adapting the proof of Proposition 2(3) of \cite{SF} to the case where the fibers of the bundles involved are annuli rather than tori).  Thus, every
orientable annulus bundle over the torus having just one boundary
component is isomorphic as a smooth fiber bundle to either $Y_0$ or
$Y_1$.

By definition, we have \begin{equation}\label{y0def} Y_0=S^1\times
\frac{\mathbb{R}\times A}{(x+1,z)\sim (x,z^{-1})}\end{equation}
Meanwhile, we see easily that \begin{equation}\label{y1}
Y_1=\frac{\mathbb{R}\times S^1\times A}{(x+1,e^{i\theta},z)\sim
(x,e^{i\theta},e^{i\theta}z^{-1})},\end{equation} since the right
hand side above obviously admits the structure of an annulus bundle
and so by our earlier remarks is isomorphic either to $Y_0$ or to
$Y_1$; computation of the fundamental group then shows that it is
distinct from $Y_0$.

\begin{lemma} \label{y0} Let $\pi\co E\to T^2$ be an $S^2$-bundle with symplectic form $\omega\in \Omega^2(E)$,
and let $\Sigma\subset E$ be an embedded, connected, symplectic
representative of $-PD(\kappa_E)\in H_2(E;\mathbb{Z})$, with tubular
neighborhood $\nu\Sigma$.  Then $E\setminus \nu\Sigma$ is
diffeomorphic to $Y_0$ if and only if $E$ is symplectomorphic to
$S^2\times T^2$ (with some split symplectic form).\end{lemma}

\begin{proof}  By Lemma \ref{annbundle}, possibly after redefining $\pi$ we may assume that $\pi$ has symplectic fibers and that $\pi|_{\Sigma\setminus \nu\Sigma}$ defines an annulus bundle over $T^2$. 

For the forward implication, simply note that the
annulus bundle $Y_0$ admits a section $(\theta,x)\mapsto
[\theta,x,1]$; if $E\setminus \nu\Sigma$ is diffeomorphic to $Y_0$
this section includes into $E$ as a torus which intersects the
fibers of $E\to T^2$ once transversely and which misses $\Sigma$.
Now the total space of $E$ is, by results of \cite{Mc1},
symplectomorphic by a fiber-preserving map to the projectivation of
a complex line bundle over $T^2$ of degree either $0$ or $1$;
however, in the projectivization of a line bundle of degree $1$ over
the torus there are no homology classes having intersection number
$1$ with the fibers and $0$ with the anticanonical class.  Hence $E$
must be the projectivization of the trivial complex line bundle over
$T^2$, \emph{i.e.} $S^2\times T^2$.

Conversely, suppose that $E$ is symplectomorphic to $S^2\times T^2$.
As in Lemma \ref{annbundle}, we can assume that $\Sigma$ meets each
fiber of $\pi\co E\to T^2$ transversely twice.  Now let $p\co M\to
S^2$ be a nontrivial $S^2$-bundle over $S^2$, and let $F\subset M$
be the disjoint union of a section of square $-1$ and a section of
square $1$ of $p$, each of which is symplectic with respect to some
symplectic form on $F$ which restricts nondegenerately to the fibers
of $\pi$.  Then by the pairwise sum construction in \cite{G}, the
fiber sum $E'$ of $M$ and $E$ carries a symplectic form and admits a
symplectic torus $\Sigma'$ obtained by gluing $\Sigma$ to the
section of square $-1$ in $F$ at one of its intersection points with
the fiber and to the section of square $1$ in $F$ at the other. Now
the induced $S^2$-fibration $\pi'\co E'\to T^2$ on the fiber sum is
easily seen to admit sections of odd square (glue a section of even
square in $E$ to a section of odd square in $M$), so $E'$ is
diffeomorphic to the nontrivial $S^2$-bundle over $T^2$.  Hence by
the previous paragraph $E'\setminus \nu\Sigma'$ is not diffeomorphic
to $Y_0$, so it is diffeomorphic to $Y_1$.

Now $M$ is diffeomorphic to
$\mathbb{C}P^2\#\overline{\mathbb{C}P^2}$; the complement
$M\setminus \nu F$ of a neighborhood of the disjoint union $F$ of a
section of square $1$ and a section of square $-1$ is then
diffeomorphic to the complement of a neighborhood of the union of a point and a line in
$\mathbb{C}P^2$, \emph{i.e.} to a region
$\{(z,w)\in\mathbb{C}P^2|r\leq |z|^2+|w|^2\leq R\}$ in
$\mathbb{C}^2$.  In these terms, $\pi|_{M\setminus \nu F}\co
M\setminus \nu F\to\mathbb{C}P^1$ is the Hopf map $(z,w)\mapsto
[z:w]$.  This shows that the annulus fibration $p|_{M\setminus \nu
F}\co M\setminus \nu F\to S^2$ is obtained from the trivial annulus
fibration over $S^2$ by removing the neighborhood of a fiber and
regluing it by the diffeomorphism $(e^{i\theta},z)\mapsto
(e^{i\theta},e^{i\theta}z)$ of $\partial D^2\times A$.

But the annulus fibration $\pi'\co E'\setminus \nu\Sigma'\to T^2$ is
obtained by taking the fiber sum of $\pi\co E\setminus \nu\Sigma\to
T^2$ with $p\co M\setminus \nu F\to S^2$, so this implies that
$\pi'\co E'\setminus\nu\Sigma'\to T^2$ may be constructed from
$\pi\co E\setminus \nu \Sigma\to T^2$ by removing the neighborhood
of a fiber and regluing it by the diffeomorphism
$(e^{i\theta},z)\mapsto (e^{i\theta},e^{i\theta}z)$ of $\partial
D^2\times A$.  Now performing this operation on the annulus bundle
$Y_0$ yields $Y_1$, while performing it on $Y_1$ yields $Y_2\cong
Y_0$. So since we've already established that $E'\setminus
\nu\Sigma'$ is diffeomorphic to $Y_1$, it must be that $E\setminus
\nu\Sigma$ is diffeomorphic to $Y_0$.\end{proof}

\begin{theorem}
Let $X$ be an $S^2$-bundle over $T^2$, let $F\subset X$ be an embedded
symplectic representative of $-PD(\kappa_X)$, and let $F'\subset
E(1)$ be an embedded symplectic representative of
$-PD(\kappa_{E(1)})$.  Then the symplectic sum $E(1)\#_{F'=F}X$ is
diffeomorphic to the Enriques surface.
\end{theorem}

\begin{proof}

First, we notice that we can reduce to the case that $X$ is
symplectomorphic to $S^2\times T^2$ (with some split symplectic
form). Indeed, if $X$ is instead diffeomorphic to the nontrivial
$S^2$-bundle over $T^2$, we shall twice apply the result of \cite{G}
and \cite{MSy} which allows us to ``trade blowups'' as discussed in
Case 1. First, if $(M_1,F_1)$ is the result of blowing down a
$(-1)$-sphere passing once positively and transversely through $F'$ (to find such a sphere, use an almost complex structure preserving $TF'$ to evaluate the Gromov--Witten invariant of one of the classes of the exceptional spheres of $E(1)=\mathbb{C}P^2\#9\overline{\mathbb{C}P^2}$)
and if $(M_2,F_2)=(X,F)$, we see
that $E(1)\#_{F'=F}X$ is deformation equivalent to
$M_{1}\#_{F_1=\tilde{F_2}}\tilde{M_2}$. Now the ruling $X\to T^2$
induces a genus-0 (not relatively minimal) Lefschetz fibration
$\pi\co \tilde{M}_2\to T^2$ each of whose fibers meets the symplectic
square-$(-1)$ torus $\tilde{F}_2$ twice; $\pi$ has just one singular fiber,
whose components $C_1$ and $C_2$ are two embedded $(-1)$-spheres
(one of which, say $C_1$, is the exceptional sphere of the blowup,
and the other of which is the proper transform of the fiber of $X\to
S^2$ that passes through the blown-up point), each of which
intersects $\tilde{F}_2$ once.  Now blowing down $C_2$ produces a manifold
symplectomorphic to $S^2\times T^2$, and $\tilde{F}_2\subset M_2$ is
isotopic to the proper transform of a symplectic representative
$F''$ of $-PD(\kappa_{S^2\times T^2})$.  Hence
$M_{1}\#_{F_1=\tilde{F_2}}\tilde{M_2}$ is in turn deformation
equivalent to $\tilde{M_1}\#_{\tilde{F_1}=F''}(S^2\times T^2)$.
Since $(\tilde{M_1},\tilde{F_1})$ is obtained by first blowing down
a sphere passing once positively and transversely through $F'\subset E(1)$ and then blowing up a
point on the image of $F'$ under the blowdown, it follows that
$\tilde{M_1}$ is deformation equivalent to $E(1)$ and $\tilde{F_1}$
represents $-PD(\kappa_{\tilde{M_1}})$.  This allows us to
hereinafter assume that $X=S^2\times T^2$.

 By Lemma \ref{isotdiff},
$E(1)\setminus \nu F'$ is diffeomorphic to the manifold with
boundary  $N$ obtained by deleting a neighborhood of a regular fiber
of an elliptic fibration on $E(1)$, while $X\setminus \nu F$ is, by
Lemma \ref{y0} and our reduction to the case that $X=S^2\times T^2$,
diffeomorphic to $Y_0$.

  So the
symplectic sum  in question is diffeomorphic to \[
X_0=N\cup_{\partial} Y_0; \] note that since by Proposition 1 of
Appendice 4 of \cite{GM} every orientation preserving diffeomorphism
of $\partial N$ extends to $N$, the diffeomorphism type of $X_0$
is determined independently of the boundary gluing maps.

We claim now that $X_0$ is diffeomorphic to the Enriques surface.
 In fact, this is essentially a remark on p. 50 of
\cite{GM}; for a direct proof, recall that $Y_0=S^1\times Z_0$ where
$Z_0=\mathbb{R}\times A/(x+1,z)\sim (x,z^{-1})$.  Now projecting
$Z_0$ onto its \emph{second} factor gives $Z_0$ the structure of a
Seifert fibration over $D^2$ with two multiple fibers each having
multiplicity $2$; hence $Y_0=S^1\times Z_0$ is the result of
performing two multiplicity two logarithmic transformations on the
trivial elliptic fibration $T^2\times D^2$. Thus
$X_0=N\cup_{\partial}Y_0$ is obtained from $E(1)$ by deleting a
neighborhood of a smooth fiber and replacing that neighborhood with
the result of two multiplicity two logarithmic transformations on
$T^2\times D^2$, \emph{i.e.}, $X_0$ is obtained from $E(1)$ by
performing two multiplicity-two logarithmic transformations.  But
this is precisely the definition of the Enriques surface.

\end{proof}

\subsection{Case 3: $X_1$ and $X_2$ are irrational and ruled}

Since a $k$-fold blowup of an $S^2$-bundle over a surface of genus
$h$ has $c_{1}^{2}=8-8h-k$, in order for the symplectic sum of
irrational ruled surfaces $X_1$ and $X_2$ along a torus to have
$c_{1}^{2}=0$, both $X_1$ and $X_2$ must be $S^2$-bundles over $T^2$. By Theorem \ref{minusf}, the surfaces $F_i$ are embedded
symplectic tori representing $-PD(\kappa_{X_i})$.  As in the proof
of Lemma \ref{annbundle}, results of \cite{Mc1} imply that there are
projections $\pi_i\co X_i\to T^2$ such that $\pi_i|_{F_i}$ is an
unramified double cover of $T^2$ by $F_i$; the deck transformation
of this cover is then a free orientation-preserving involution $\tau_i\co F_i\to F_i$.  By
considering these involutions, we shall realize any symplectic sum
of the $X_i$ along the $F_i$ as the total space of some torus bundle
over $T^2$.

\begin{lemma}\label{isotinv} Let $\tau_1,\tau_2\co T^2\to T^2$ be free orientation-preserving involutions,
and let $\phi\co T^2\to T^2$ be any diffeomorphism.  Then $\phi$ is
isotopic to a diffeomorphism $\phi'\co T^2\to T^2$ with the property
that
\[ \phi'^{-1}\circ \tau_2\circ\phi'\circ\tau_1 \mbox{ is either the
identity or a free involution.}\]

\end{lemma}

\begin{proof} First of all, note that any two free orientation-preserving involutions
$\tau,\tau'\co T^2\to T^2$ are conjugate.  Indeed, letting $E$ be
the quotient of $T^2$ by $\tau$, $E'$ the quotient of $T^2$ by
$\tau'$, and $\pi\co T^2\to E$, $\pi'\co T^2\to E'$ the projections,
$E$ and $E'$ are both tori, so that there exists a diffeomorphism
$\psi\co E\to E'$.  The images of $\pi_1(T^2)$ in $\pi_1(E')$ by
$\pi'$ and $\psi\circ \pi$ are both index $2$ lattices in
$\pi_1(E')\cong \mathbb{Z}^2$, so there is an element $A$ of
$SL(2;\mathbb{Z})$ taking one to the other; hence by composing
$\psi$ with a diffeomorphism of $E'$ that induces $A$ on $\pi_1$ we
can assume that the maps induced on $\pi_1$ by $\pi'$ and by
$\psi\circ\pi$ have the same image.  Hence $\psi\circ \pi\co T^2\to
E'$ lifts to a diffeomorphism $f\co T^2\to T^2$ such that $\pi'\circ
f=\psi\circ\pi$.  Since $\tau$ (resp. $\tau'$) takes $x\in T^2$ to
the unique other point in $\pi^{-1}(\pi(x))$ (resp.
$\pi'^{-1}(\pi'(x))$) it follows that $\tau'\circ f=f\circ\tau$, so
$\tau'$ and $\tau$ are indeed conjugate.

In light of this, identifying $T^2=\mathbb{R}^2/\mathbb{Z}^2$ and
conjugating $\tau_1,\tau_2,\phi$ by some diffeomorphism, we can
assume that $\tau_1([x,y])=[x+1/2,y]$ (where $[x,y]\in T^2$ is the
equivalence class of $(x,y)\in\mathbb{R}^2$ under the relations
$(x+1,y)\sim (x,y+1)\sim (x,y)$).  By the previous paragraph, since
$\phi^{-1}\circ\tau_2\circ\phi$ and $\tau_1$ are free involutions,
there is some $\alpha\in Diff(T^2)$ such that
$\phi^{-1}\circ\tau_2\circ\phi=\alpha^{-1}\circ\tau_1\circ\alpha$.
Now $\alpha$ is isotopic to some linear diffeomorphism
$A=\left(\begin{array}{cc}a&b\\c&d\end{array}\right)\in
SL(2;\mathbb{Z})$; say $A=\alpha\circ f_1$ where $\{f_t\}_{t\in
[0,1]}$ is a smooth family of diffeomorphisms such that $f_0=1$.
Then where $\phi'=\phi\circ f_1$, $\phi$ is isotopic to $\phi'$
and we have $\phi'^{-1}\circ\tau_2\circ\phi'=A^{-1}\circ\tau_1\circ
A$.

Now since $\tau_1([x,y])=[x+1/2,y]$, one easily computes \[
A^{-1}\circ\tau_1\circ
A\circ\tau_1([x,y])=\big[x+\frac{d+1}{2},y-\frac{c}{2}\big],\] which
defines the identity if  $c$ is even (forcing $d$ to be odd since
$A\in SL(2;\mathbb{Z})$) and a free involution if $c$ is odd.
\end{proof}

\begin{theorem} \label{makebundle} Let $\pi_i\co X_i\to T^2$ $(i=1,2)$ be
$S^2$-bundles over $T^2$, $F_i\subset X_i$ embedded tori with the
property that $\pi_i|_{F_i}$ is an unramified double cover.  Let
$\nu_i$ be tubular neighborhoods of $F_i$ (each identified with
$D^2\times F_i$), and let $\Phi\co
\partial\nu_1\to\partial\nu_2$ be a diffeomorphism which
(viewing $\partial\nu_i$ as an $S^1$-bundle over $F_i$) covers
some diffeomorphism $\phi\co F_1\to F_2$. Then the normal connect
sum \[ (X_1\setminus \nu_1)\cup_{\partial\nu_1\sim_{\Phi}\partial\nu_2}(X_2\setminus \nu_2)\] is
diffeomorphic to the total space of a $T^2$-bundle over $T^2$.
\end{theorem}

\begin{proof}
First, note that (after performing isotopies which don't change the
diffeomorphism type of the normal connect sum), we can assume that,
for $i=1,2$, the $S^2$-bundle projection $\pi_i$ is constant on each
fiber of the disc bundle projection $\nu_i\to F_i$, and that
(using Lemma \ref{isotinv}) $\phi^{-1}\circ
\tau_2\circ\phi\circ\tau_1$ is either the identity or a free
involution, where $\tau_i\co F_i\to F_i$ is the deck transformation
induced by the cover $\pi_i|_{F_i}$. Let $Z=(X_1\setminus \nu_1)\cup_{\partial\nu_1\sim_{\Phi}\partial\nu_2}(X_2\setminus
\nu_2)$.

Suppose that $\phi^{-1}\circ \tau_2\circ\phi\circ\tau_1$ is the
identity.  We define a bundle map $\pi\co Z\to T^2$ as follows.  If
$x\in X_1\setminus \nu_1\subset Z$, put $\pi(x)=\pi_1(x)$.  If
$x\in X_2\setminus\nu_2\subset Z$, then there are two points
$x_2,\tau_2(x_2)\in F_2\cap\pi_{2}^{-1}(\{\pi_2(x)\})$, and since
$\phi^{-1}\circ \tau_2\circ\phi\circ\tau_1$ is the identity we have
$\tau_1(\phi^{-1}(x_2))=\phi^{-1}(\tau_2(x_2))$, so that
$\pi_1(\phi^{-1}(x_2))=\pi_1(\phi^{-1}(\tau_2(x_2)))$ and we set
\[ \pi(x)=\pi_1(\phi^{-1}(x_2))=\pi_1(\phi^{-1}(\tau_2(x_2))).\]
Since, for each $p\in F_i$, the fiber of the circle bundle
$\partial\nu_i\to F_i$ over $p$ is mapped by $\pi_i$ to $\pi_i(p)$, our map
$\pi$ is defined consistently on the identified boundary components
$\partial\nu_1,\partial\nu_2$ in $Z=(X_1\setminus \nu_1)\cup_{\partial\nu_1\sim_{\Phi}\partial\nu_2}(X_2\setminus
\nu_2)$.  One easily sees that $\pi\co Z\to T^2$ is a
$T^2$-fibration; the point here is that since $\phi^{-1}\circ
\tau_2\circ\phi\circ\tau_1=1$, $\phi\co F_1\to F_2$ descends to a
map $f\co T^2\to T^2$ such that $\pi_2|_{F_2}\circ
\phi=f\circ\pi_1|_{F_1}$; the fiber of $\pi$ over $t\in T^2$ is
formed by gluing the annulus $\pi_{1}^{-1}(\{t\})\cap (X_1\setminus
\nu_1)$ to the annulus $\pi_{2}^{-1}(\{f(t)\})\cap (X_2\setminus
\nu_2)$.

It remains to consider the case that $\phi^{-1}\circ
\tau_2\circ\phi\circ\tau_1$ is a free involution.  Then
$\phi^{-1}\circ\tau_2\circ\phi\circ\tau_1$ commutes with $\tau_1$
and their composition (namely $\phi^{-1}\circ\tau_2\circ\phi$) is
also a free involution.  Let $E=F_1/\langle \tau_1,\phi^{-1}\circ
\tau_2\circ\phi\circ\tau_1\rangle$ and let $p\co F_1\to E$ be the
projection (which is an unramified quadruple covering of a torus by
a torus). We shall define a torus fibration $\pi\co Z\to E$ rather
similarly to the previous case, except that here the fibers will be
formed by gluing four annuli rather than two.  If $x\in X_1\setminus
\nu_1\subset Z$, set $\pi(x)=p(x_1)$ where $x_1\in F_1\cap
\pi_{1}^{-1}(\{\pi_1(x)\})$; this is a coherent definition since the
two elements of $\pi_{1}^{-1}(\{\pi_1(x)\})$ are intertwined by
$\tau_1$. If $x\in X_2\setminus \nu_2\subset Z$, we intend to set
$\pi(x)=p(\phi^{-1}(x_2))$ where $x_2\in
F_2\cap\pi_2^{-1}(\{\pi_2(x)\})$; we need to see that the two
possible choices of $x_2$ (either of which is taken to the other by
$\tau_2$) give the same value for $p(x)$.  In other words, we need
to check that if $x_2\in F_2$ then
$p(\phi^{-1}(x_2)))=p(\phi^{-1}(\tau_2(x_2)))$.  Now since
$\phi^{-1}\circ \tau_2\circ\phi\circ\tau_1$ has order $2$,
\[
\phi^{-1}(\tau_2(\phi(\tau_1(\phi^{-1}(x_2)))))=\tau_1(\phi^{-1}(\tau_2(x_2)))
\] so since $E$ is the quotient of $F_1$ by $\tau_1$ and
$\phi^{-1}\circ\tau_2\circ\phi\circ\tau_1$ it is indeed the case
that $p(\phi^{-1}(x_2)))=p(\phi^{-1}(\tau_2(x_2)))$ for each $x_2\in
F_2$.  We have thus defined $\pi\co Z\to E$; it is again easily seen
to be a torus bundle, with its fibers of the shape \[
\frac{A_0\coprod A_1\coprod A_2\coprod
A_3}{\partial_+A_i\sim\partial_-A_{i+1} (i\in
\mathbb{Z}/4\mathbb{Z})}
\] where $A_0$ and $A_2$ are annulus fibers of
$\pi_2|_{X_2\setminus\nu_2}$ and $A_1$ and $A_3$ are annulus
fibers of $\pi_1|_{X_1\setminus\nu_1}$ (the fact that
$\phi^{-1}\circ\tau_2\circ\phi\circ\tau_1$ is free serves to ensure
that, in each of these torus fibers, $A_0$ and $A_2$ are distinct,
as are $A_1$ and $A_3$).
\end{proof}

This shows that any symplectic $4$-manifold obtained as the
symplectic sum of two $S^2$-bundles over $T^2$ along a pair of
bi-sections is diffeomorphic to a $T^2$-bundle over $T^2$.  In fact,
we can be quite specific about which $T^2$-bundles over $T^2$ are
obtained in this fashion.  $T^2$-bundles over $T^2$ were classified
in \cite{SF}; in particular, Theorem 5 of that paper shows that the
total spaces of such bundles are distinguished from one another up
to diffeomorphism by their fundamental groups.  As such, finding the
diffeomorphism type of the manifold $(X_1\setminus \nu_1)\cup_{\partial\nu_1\sim_{\Phi}\partial\nu_2}(X_2\setminus
\nu_2)$ in Theorem \ref{makebundle} is just a matter of applying
van Kampen's theorem.

We know that, for $i=1,2$, the manifold $X_i\setminus \nu_i$ is
diffeomorphic to one of the manifolds $Y_0$ or $Y_1$ of
(\ref{y0def}),(\ref{y1}); more specifically, if $X_i$ is
diffeomorphic to $S^2\times T^2$ then $X_i\setminus \nu_i\cong
Y_0$, and otherwise $X_i\setminus \nu_i\cong Y_1$.  Note that
\begin{equation}\label{annpi1}
\pi_1(Y_j)=\frac{\langle\alpha,\beta,m\rangle}{\alpha^{-1}m\alpha=m^{-1},\beta
m=m\beta,\alpha\beta\alpha^{-1}\beta^{-1}=m^j}\,
(j=0,1),\end{equation} with $\partial Y_j$ being spanned by the
subgroup generated by $\alpha^2,\beta,m$.  $m$ here is the generator
of the fundamental group of the annulus fiber of the bundle map
$Y_j\to T^2$. Where $F$ denotes the torus $F_1$ or $F_2$ whose neighborhood we have removed from $X_1$ or $X_2$ to get $Y_j$, we have a trivial circle bundle $p_j\co \partial Y_j\to F$ whose action on $\pi_1$ has kernel $\langle m\rangle$.   There is, of course, some flexibility in the choice
of the generators: first, we get the same presentation if we replace $\beta$ by $\beta m^j$ and then $m$ by $m^{-1}$; secondly, if $\left(\begin{array}{cc}p&q\\r&s\end{array}\right)\in SL(2,\mathbb{Z})$ and $q$ is even (so that $p$ and $s$ are odd; say $s=2t+1$), then we get the same presentation by replacing $\alpha$ with $\alpha'=\alpha^p\beta^r$ and $\beta$ with $\beta'=\alpha^q\beta^sm^{jt}$.
% , corresponding to changing the basis for $\pi_1(F)$: in particular, the generators $\alpha,\beta$
%could be replaced in the presentation (\ref{annpi1}) by any elements
%$\alpha',\beta'$ with the property that $\langle
%\alpha'^2,\beta'\rangle=\langle\alpha^2,\beta\rangle\leq
%\pi_1(Y_j)$, that is to say, by any pair
%$\alpha'=\alpha^p\beta^r,\beta'=\alpha^q\beta^s$ where \[
%\left(\begin{array}{cc}p&q\\r&s\end{array}\right)\in
%SL(2;\mathbb{Z}) \mbox{ and } q\mbox{ is even.}\]
One convenient
consequence of this is that if $\gamma\in \pi_1(\partial Y_j)$ is any element with the property that $(p_j)_*\gamma$
is primitive in $\mathbb{Z}^2=\pi_1(F)$ then  the generators $\alpha,\beta$ in the presentation
(\ref{annpi1}) may be chosen so that $\gamma$ takes one of the forms
\[ \gamma=\alpha^2\beta^{2c}m^e\mbox{ or } \gamma=\alpha^{2a}\beta
m^e.\]

We now consider the manifold resulting from gluing two of these
manifolds $Y_j,Y_k$ ($j,k\in\{0,1\}$) together along their
boundaries in a way consistent with the symplectic sum operation.
Now in terms of bases $\{\alpha_{1}^2,\beta_1,m_1\}$,
$\{\alpha_{2}^{2},\beta_2,m_2\}$ for the fundamental groups of the
boundaries $\partial Y_j$ and $\partial Y_k$ respectively, since the
gluing map $\Phi$ is required to cover an isomorphism of the normal
bundles it will identify $m_1$ with $m_{2}$ (possibly after replacing one of the $m_i$ with its inverse, which as mentioned earlier can be done without affecting the presentation (\ref{annpi1}) at the cost of possibly multiplying $\beta_i$ by $m_i$); also
$\alpha_{1}^{2}$, since it projects via $p_j$ to a primitive element in $\pi_1(F)$,
will be taken to some element in the fundamental group of $\partial Y_k$ which likewise projects via $p_k$ to a primitive element. Hence by the remark at the end
of the previous paragraph, possibly after renaming the generators
$\alpha_2$, $\beta_2$, and $m_2$ in the presentation of $\pi_1(Y_k)$, the
action of the gluing map on the fundamental groups of the boundaries in terms
of the bases $\{\alpha_{i}^{2},\beta_i,m_i\}$  takes one of the
forms
\[ \left(\begin{array}{ccc}1&b&0\\2c&d&0\\e&f&1\end{array}\right)
\,(d-2bc=1)\mbox{ or }
\left(\begin{array}{ccc}a&b&0\\1&d&0\\e&f&1\end{array}\right)\,
(ad-b=1).\]  Hence  van Kampen's theorem gives the fundamental group
of the glued manifold $Y_j\cup_{\Phi}Y_k$ as either
\begin{equation} \label{c-even} \pi_1(Y_j\cup_{\Phi}Y_k)=\frac{\langle
\alpha_1,\beta_1,\alpha_2,\beta_2,m\rangle}{\begin{array}{c}\alpha_{1}^{-1}m\alpha_1=\alpha_{2}^{-1}m\alpha_2=m^{-1},\,\beta_1m=m\beta_1,\beta_2m=m\beta_2,\alpha_1\beta_1\alpha_{1}^{-1}\beta_{1}^{-1}=m^j\\
 \alpha_2\beta_2\alpha_{2}^{-1}\beta_{2}^{-1}=m^k,
\,\alpha_{1}^{2}=\alpha_{2}^{2}\beta_{2}^{2c}m^e,\,\beta_1=\alpha_{2}^{2b}\beta_{2}^{d}m^f\end{array}}\end{equation}
or
\begin{equation}\label{c-odd}\pi_1(Y_j\cup_{\Phi}Y_k)=\frac{\langle
\alpha_1,\beta_1,\alpha_2,\beta_2,m\rangle}{\begin{array}{c}\alpha_{1}^{-1}m\alpha_1=\alpha_{2}^{-1}m\alpha_2=m^{-1},\,\beta_1m=m\beta_1,\beta_2m=m\beta_2,\,\alpha_1\beta_1\alpha_{1}^{-1}\beta_{1}^{-1}=m^j
\\ \alpha_2\beta_2\alpha_{2}^{-1}\beta_{2}^{-1}=m^k,
\,\alpha_{1}^{2}=\alpha_{2}^{2a}\beta_{2}m^e,\,\beta_1=\alpha_{2}^{2b}\beta_{2}^{d}m^f\end{array}}\end{equation}

The reader may verify that the group on the right hand side of
(\ref{c-even}) may be rewritten, by identifying
$\gamma=\alpha_{1}^{-1}\alpha_2\beta_{2}^{c}$, as \[ \frac{\langle
\alpha_1,\beta_2,m,\gamma\rangle}{\begin{array}{c}m\gamma=\gamma
m,\,\alpha_1\beta_2\alpha_{1}^{-1}\beta_{2}^{-1}=m^{j+2(f-be)},\,\alpha_{1}^{-1}m\alpha_1=m^{-1},\\
\alpha_{1}^{-1}\gamma\alpha_1=m^{kc-e}\gamma^{-1},\,\beta_{2}^{-1}m\beta_2=m,\,\beta_{2}^{-1}\gamma\beta_2=m^{j-k+2(f-be)}\gamma\end{array}},\]
which we recognize as the fundamental group of the $T^2$-bundle over
$T^2$ given in the notation of the introduction as \[
M\left(\left(\begin{array}{cc}-1&kc-e\\0&-1\end{array}\right),\left(\begin{array}{cc}1&j-k+2(f-be)\\0&1\end{array}\right);\left(\begin{array}{c}j+2(f-be)\\0\end{array}\right)\right).\]

Similarly, the group on the right hand side of (\ref{c-odd}) may be
identified, by taking
$\sigma=\alpha_1\alpha_2\alpha_{1}^{-1}\alpha_{2}^{-1}$ and then
using the relations $\alpha_1\alpha_{2}^{2}=\alpha_{2}^{2}\alpha_1
m^{j+2(f-de)}$ and $\alpha_2\alpha_{1}^{2}=\alpha_{1}^{2}\alpha_2
m^{2e-k}$ to obtain commutation relations between the $\alpha_i$ and
$\sigma$, as \[ \frac{\langle
\alpha_1,\alpha_2,m,\sigma\rangle}{\begin{array}{c}m\sigma=\sigma
m,\,\alpha_1\alpha_2\alpha_{1}^{-1}\alpha_{2}^{-1}=\sigma,\,\alpha_{1}^{-1}m\alpha_1=m^{-1},\\
\alpha_{1}^{-1}\sigma\alpha_1=m^{k-2e}\sigma^{-1},\,\alpha_{2}^{-1}m\alpha_2=m^{-1},\,\alpha_{2}^{-1}\sigma\alpha_2=m^{j+2(f-de)}\sigma^{-1}\end{array}},\]
which is precisely the fundamental group of the $T^2$-bundle over
$T^2$ \[
M\left(\left(\begin{array}{cc}-1&k-2e\\0&-1\end{array}\right),\left(\begin{array}{cc}-1&j+2(f-de)\\0&-1\end{array}\right);\left(\begin{array}{c}0\\1\end{array}\right)\right).\]

Now by changing the basis for the homology of the base by
$\left(\begin{array}{cc}p&q\\r&s\end{array}\right)\in
SL(2;\mathbb{Z})$, a $T^2$-bundle over $T^2$ of form
$M(A,B;\vec{v})$ may be equated with $M(A^pB^r,A^qB^s;\vec{v})$;
also, the bundles $M(A,B;\vec{v})$ and $M(A,B;\vec{v}')$ are
equivalent if $\vec{v}'-\vec{v}$ lies in the submodule of
$\mathbb{Z}^2$ spanned by the columns of $A-I$ and $B-I$ (where $I$
is the identity; these statements are proven in Proposition 2 of
\cite{SF}).  As such, given a bundle of form \[
M\left(\left(\begin{array}{cc}-1&\delta\\0&-1\end{array}\right),\left(\begin{array}{cc}1&\zeta\\0&1\end{array}\right);\left(\begin{array}{c}j+2x\\0\end{array}\right)\right),\]
by letting $z=gcd(\delta,\zeta)$, $p=\zeta/z$, $r=\delta/z$, and
(since $p$ and $r$ are then relatively prime) $q$ and $s$ be such
that $ps-qr=1$, so that $-q\delta+s\zeta=z$, we obtain \[
M\left(\left(\begin{array}{cc}-1&\delta\\0&-1\end{array}\right),\left(\begin{array}{cc}1&\zeta\\0&1\end{array}\right);\left(\begin{array}{c}j+2x\\0\end{array}\right)\right)=
M\left((-1)^pI,(-1)^q\left(\begin{array}{cc}1&z\\0&1\end{array}\right);\left(\begin{array}{c}j\\0\end{array}\right)\right)\]
(note also that $p$ and $q$ can't both be even since $ps-qr=1$, and
if $p$ is odd then a further basis change for the homology of the
base identifies $M(-I,-A;\vec{v})$ with $M(-I,A;\vec{v})$).  This
gives rise to the following list of possibilities for the
diffeomorphism type of $Y_j\cup_{\Phi}Y_k$ when its fundamental
group is given by (\ref{c-even}):
\begin{center}
\begin{tabular}{l|l|l} $j$&$k$&possible diffeomorphism types of
$Y_j\cup_{\Phi}Y_k$ \\\hline $0$& $0$ &
$M\left(I,\left(\begin{array}{cc}-1&z\\0&-1\end{array}\right);\left(\begin{array}{c}0\\0\end{array}\right)\right),\,M\left(-I,\left(\begin{array}{cc}1&2y\\0&1\end{array}\right);\left(\begin{array}{c}0\\0\end{array}\right)\right)$
$(y,z\in\mathbb{Z})$\\ \hline $0$& $1$ &
$M\left(-I,\left(\begin{array}{cc}1&2y+1\\0&1\end{array}\right);\left(\begin{array}{c}0\\0\end{array}\right)\right)$
$(y\in\mathbb{Z})$\\ \hline $1$& $1$ &
$M\left(I,\left(\begin{array}{cc}-1&z\\0&-1\end{array}\right);\left(\begin{array}{c}1\\0\end{array}\right)\right),\,M\left(-I,\left(\begin{array}{cc}1&2y\\0&1\end{array}\right);\left(\begin{array}{c}1\\0\end{array}\right)\right)$
$(y,z\in\mathbb{Z})$
\end{tabular}
\end{center}
Similarly, a bundle of form \[
M\left(\left(\begin{array}{cc}-1&\delta\\0&-1\end{array}\right),\left(\begin{array}{cc}-1&\zeta\\0&-1\end{array}\right);\left(\begin{array}{c}0\\1\end{array}\right)\right)\]
is equivalent to \[
M\left((-1)^{p+r}I,(-1)^{q+s}\left(\begin{array}{cc}1&
z\\0&1\end{array}\right);\left(\begin{array}{c}0\\1\end{array}\right)\right)
\] (where $z=gcd(\delta,\zeta)$, $p=\zeta/z$, $r=-\delta/z$, and
$ps-qr=1$).  From this, we deduce the following list of
possibilities for the diffeomorphism type of $Y_j\cup_{\Phi}Y_k$
when its fundamental group is given by (\ref{c-odd}):
\begin{center}
\begin{tabular}{l|l|l} $j$&$k$&possible diffeomorphism types of
$Y_j\cup_{\Phi}Y_k$ \\\hline $0$& $0$ &
$M\left(I,\left(\begin{array}{cc}-1&2y\\0&-1\end{array}\right);\left(\begin{array}{c}0\\1\end{array}\right)\right),\,M\left(-I,\left(\begin{array}{cc}1&2y\\0&1\end{array}\right);\left(\begin{array}{c}0\\1\end{array}\right)\right)$
$(y\in\mathbb{Z})$\\ \hline $0$& $1$ &
$M\left(-I,\left(\begin{array}{cc}1&2y+1\\0&1\end{array}\right);\left(\begin{array}{c}0\\1\end{array}\right)\right)$
$(y\in\mathbb{Z})$\\ \hline $1$& $1$ &
$M\left(I,\left(\begin{array}{cc}-1&2y+1\\0&-1\end{array}\right);\left(\begin{array}{c}0\\1\end{array}\right)\right)$
$(y\in\mathbb{Z})$
\end{tabular}
\end{center}
In both of the above tables, it's easy to see that any of the
indicated diffeomorphism types can in fact be realized by means of
an appropriate choice of the gluing map $\Phi$.  Since if $X$ is a
ruled surface over $T^2$ and $F\subset X$ is an embedded
representative of $-PD(\kappa_X)$, we've seen that $X\setminus \nu
F\cong Y_0$ if $X=S^2\times T^2$ and $X\setminus \nu F\cong Y_1$ if
$X=S^2\tilde{\times}T^2$, this completes the proof that the
diffeomorphism types of the symplectic sums in question are as
claimed in the introduction.

\end{document}